\documentclass[10pt,reqno]{amsart}
\usepackage{amsfonts,amsthm,latexsym,amsmath,amssymb,amscd,amsmath}
\usepackage{graphics,graphicx}
\usepackage{epic}
\usepackage{eepic}
\usepackage{xypic}

\theoremstyle{plain}
   \newtheorem{theorem}{Theorem}[section]
   \newtheorem{proposition}[theorem]{Proposition}
   \newtheorem{lemma}[theorem]{Lemma}
   \newtheorem{corollary}[theorem]{Corollary}
   \newtheorem{problem}{Problem}[section]
   \newtheorem{conjecture}[problem]{Conjecture}
   \newtheorem{question}[problem]{Question}
\theoremstyle{definition}
   \newtheorem{definition}{Definition}[section]
   \newtheorem{example}{Example}[section] 
   
\theoremstyle{remark}
   \newtheorem{remark}{Remark}[section]

\newcommand{\LC}{\text{LC}}
\newcommand{\SLC}{\text{SLC}}
\newcommand{\NN}{\mathbb{N}}
\newcommand{\PP}{\mathcal{P}}
\newcommand{\cI}{\mathcal{I}}
\newcommand{\cA}{\mathcal{A}}
\newcommand{\AS}{\mathcal{A}}
\newcommand{\SR}{\mathcal{SR}}
\newcommand{\BS}{\mathcal{B}}

\newcommand{\PHR}{\text{PHR}}

\newcommand{\CNA}{\text{CNA}}
\newcommand{\NLC}{\text{NLC}}
\newcommand{\hNLC}{\text{h-NLC}}
\newcommand{\NA}{\text{NA}}
\newcommand{\ULC}{\text{ULC}}
\newcommand{\TP}{\text{TP}}
\newcommand{\GKK}{\text{GKK}}
\newcommand{\HH}{\mathcal{H}}

\newcommand{\RR}{\mathbb{R}}
\newcommand{\CC}{\mathbb{C}}

\renewcommand{\Im}{\mathfrak{Im}}

\newcommand{\al}{\alpha}

\newcommand{\De}{\Delta}

\newcommand{\cN}{{[n]}}

\newcommand{\ben}{\mathbf{1}}

\newcommand{\fP}{\mathfrak{P}}
\newcommand{\pa}{\partial}
\newcommand{\suc}{\succcurlyeq}
\newcommand{\pc}{\preccurlyeq}

\newcommand{\fS}{\mathfrak{S}}

\def\newop#1{\expandafter\def\csname #1\endcsname{\mathop{\rm
#1}\nolimits}}

\newop{diag}
\newop{supp}
\newop{JP}
\newop{Sym}
\newop{si}
\newop{glb}
\newop{Trp}
\newop{roots}
\newop{edges}

\begin{document}

\title[Negative Dependence and the Geometry of Polynomials]{Negative 
Dependence and \\ the Geometry of Polynomials}

\author[J.~Borcea]{Julius Borcea}
\address{Department of Mathematics, Stockholm University, SE-106 91 Stockholm,
Sweden}
\email{julius@math.su.se}

\author[P.~Br\"and\'en]{Petter Br\"and\'en}
\address{Department of Mathematics, Royal Institute of Technology, 
SE-100 44 Stockholm, Sweden}
\email{pbranden@math.kth.se}

\author[T.~M.~Liggett]{Thomas M.~Liggett}
\address{Department of Mathematics, University of California, Los Angeles, 
CA 90095-1555, USA}
\email{tml@math.ucla.edu}


\begin{abstract}
We introduce the class of {\em strongly Rayleigh} probability measures 
by means of geometric properties of their generating
polynomials that amount to the stability of the latter. This
class covers important models such as determinantal measures 
(e.g.~product measures, 
uniform random spanning tree measures) and 
distributions for symmetric exclusion processes. 
We show that strongly Rayleigh measures enjoy all virtues of negative 
dependence and we also prove a series of conjectures due to Liggett, Pemantle,
and Wagner, respectively. Moreover, we extend Lyons' recent results on 
determinantal measures and we construct counterexamples to
several conjectures of Pemantle and Wagner on negative dependence and 
ultra log-concave rank sequences.  
\end{abstract}

\maketitle

\tableofcontents

\section{Introduction}\label{s1}

Let $\mu : 2^{[n]} \rightarrow \RR$, $[n]=\{1,\ldots,n\}$, be a function attaining nonnegative values and satisfying $\sum_{S \in 2^{[n]}}\mu(S)=1$. The function $\mu$ is said to satisfy the 
{\em positive lattice condition} \cite{FKG,P} if 
\begin{equation*}
\mu(S)\mu(T) \leq \mu(S\cup T)\mu(S\cap T)\tag*{(PLC)}
\end{equation*}
for all $S,T \subseteq [n]$. The corresponding probability measure (also denoted by $\mu$) on $2^{[n]}$ defined by $\mu(\mathcal{A})=\sum_{S \in \mathcal{A}}\mu(S)$ is said to be {\em positively associated}
if
$$\int Fd\mu \int Gd\mu \le \int FG d\mu$$
for any pair of increasing functions $F,G$ on $2^{[n]}$. The latter is a
strong correlation inequality that yields many others as well as distributional
limit theorems for models of statistical mechanics such as the ferromagnetic
Ising model or certain urn models \cite{New}. A fundamental result in the
theory of positively associated random variables, 
namely the FKG theorem \cite{FKG},
asserts that PLC implies positive association. This is a powerful tool that
allows verification of (global) correlation inequalities from the PLC
property, which is a local condition and therefore often easier to check.  

There are many important examples of 
{\em negatively dependent} ``repelling'' random
variables in probability theory, combinatorics, stochastic processes and
statistical mechanics: uniform random spanning tree measures \cite{BP}, 
symmetric exclusion processes \cite{L0,L2}, random cluster models 
(with $q<1$) \cite{grimmett,Hag,P},
balanced and Rayleigh matroids \cite{CW,feder,semple-welsh,sok,W}, 
competing urns models \cite{DR}, etc; see, e.g., \cite{J-DP,P} and references 
therein for a discussion of some of these examples and several others. 
To add to this list, in \S \ref{s3-ex} we show that both the inequalities 
characterizing multi-affine real 
stable polynomials \cite{BBS1,BBS2,PB} and Hadamard-Fischer-Kotelyansky type
inequalities in matrix theory \cite{FJ,HS,HM} may in fact be viewed as natural
manifestations of negative dependence properties.

The reverse of the PLC inequality gives the so-called
{\em negative lattice condition} 
\begin{equation*}
\mu(S)\mu(T) \geq \mu(S\cup T)\mu(S\cap T)\tag*{(NLC)}
\end{equation*}
for all $S,T \subseteq [n]$, while the usual 
definition of {\em negative association} for a probability 
measure $\mu$ on $2^{[n]}$ is the ``negative'' analog of positive 
association \cite{P}, namely the inequality
$$
\int Fd\mu \int Gd\mu \ge \int FG d\mu
$$
for any pair of increasing functions $F,G$ on $2^{[n]}$ depending on 
{\em disjoint} sets of coordinates (note that the latter condition is 
important since by the Cauchy-Schwarz inequality $F$ and $F$ are always
positively correlated).
However, negative 
dependence is not nearly as robust as positive dependence. In particular,
the NLC inequality 
{\em fails} to imply negative association 
(see Example \ref{ex1}), and there is
no known local-to-global tool comparable to the FKG theorem in the
theory of positive association. Thus, no corresponding theory
of negatively dependent events exists as yet but, as explained in \cite{P},
there certainly is a need for such a theory. 
In {\em op.~cit.~}Pemantle made a systematic
study of measures on Boolean lattices with various negative dependence
properties. He proved certain relations between them and conjectured numerous 
others, so as to pave the way for a theory of negative dependence. 
However, most open problems and conjectures are unresolved to this day, and 
the scope of the positive results so far is quite restricted.

In this paper we introduce the class of {\em strongly Rayleigh} probability 
measures by means of purely geometric conditions on the 
zero sets of their generating polynomials that amount to the real
stability of the latter (Definitions \ref{d-rs-pol}--\ref{def-rsm}). 
As we show in \S \ref{s3-ex}, this class 
contains several important examples such as 
distributions at time $t\ge 0$ for symmetric exclusion
processes generated by product or deterministic measures \cite{L2} 
(or, more generally, by strongly Rayleigh measures), 
and determinantal probability measures induced by positive contractions
\cite{BKPV,Lyons}. 
Measures of the latter type include e.g.~uniform random spanning tree 
measures \cite{BP,Lyons} and arise naturally in a variety of contexts 
pertaining to fermionic/determinantal point processes and their continuous 
scaling limits, eigenvalues of random matrices, number theory,
non-intersecting paths, transfer current matrices, orthogonal polynomial
ensembles, etc; see the discussion in \S \ref{ss-det-m} and 
references therein. 
Note that the topological dimension
of the class of strongly Rayleigh probability measures on $2^{[n]}$ is $2^n-1$ 
(\cite[Lemma~4.1]{BB}) whereas 
the class of determinantal probability measures has 
dimension $n^2$, so the 
former is much larger than the latter. Of course, there are also  
many examples of measures that satisfy some negative dependence property but 
are not in the strongly Rayleigh class, e.g.~Counterexample 1 in 
\S \ref{s-countex} and the random cluster measure 
\cite{grimmett} whenever the graph contains a cycle, see \S \ref{ss-trees}. 

The geometry underlying strongly Rayleigh measures provides us with 
powerful tools -- 
such as fundamental complex analytic and geometric
results of G\aa rding and Grace-Szeg\"o-Walsh, combinatorial/probabilistic
results of Feder-Mihail, matrix theoretic methods, and recent developments in 
the theory of
stable polynomials \cite{BBS,BBS1,BBS2,PB} -- that we use in 
\S \ref{s-stable} to develop a complete theory of negative dependence for
this class of measures. Indeed, we show that 
strongly Rayleigh measures enjoy all virtues of negative
dependence, including the strongest form of negative association (CNA+). 
In particular, this allows us to prove 
several conjectures made by Liggett \cite{L2}, Pemantle \cite{P}, and 
Wagner \cite{W}, 
respectively, and to recover and extend Lyons' main results \cite{Lyons} on 
negative
association and stochastic domination for determinantal 
probability measures induced by positive contractions. Moreover, we define a
partial order on the set of strongly Rayleigh measures (by means of the 
notion of proper position for multivariate stable polynomials studied in
\cite{BBS,BBS1,BBS2,PB}), and use it to settle Pemantle's questions and
conjectures on stochastic domination for truncations of
``negatively dependent'' measures \cite{P}. 

A series of conjectures, first appearing in print in \cite[Conjecture
4]{P} but of obscure folklore origin, states that the various negative
dependent properties studied in {\em loc.~cit.~}give rise to so-called
ultra log-concave rank sequences (see \S \ref{s21} below for the definition). 
Subsequently, in \cite{W} one of these conjectures was coined ``The Big
Conjecture'' since it would imply the validity of Mason's long-standing
conjecture in enumerative graph/matroid theory \cite{Mason} for a large 
class of matroids.
In \S \ref{s-countex} we construct the first counterexamples in the literature 
to all of these conjectures except Mason's, which once again confirms the 
delicate nature of negative dependence. After this work was made publically 
available on \texttt{www.arxiv.org} other counterexamples to some 
of these conjectures have been reported \cite{KN,Mark}. 
Markstr\"om \cite{Mark} proved that negative association does not imply 
unimodality and Kahn-Neiman \cite{KN} later showed that not even strong 
negative association (CNA+) implies unimodality.  
Note that the rank sequence of the CNA+ measure constructed in 
Counterexample 1 of \S \ref{s-countex} is neither ultra nor strong 
log-concave but it is log-concave.

The problem of describing natural
negative dependence properties that are preserved by symmetric exclusion 
evolutions has 
attracted some attention in the theory of 
interacting particle systems
and Markov processes \cite{L0,L00,L2}. In \S \ref{n-s5} we prove that the 
class of strongly Rayleigh measures is invariant
under partial symmetrization (this and several other properties fail
for e.g.~Rayleigh measures). As a consequence, in \S \ref{n-s6} we provide 
an answer to the aforementioned problem and show that if the
initial distribution of a symmetric exclusion process is strongly Rayleigh, 
then so is the distribution at time $t \geq 0$; therefore,  
by the results in \S \ref{s-sna}, the latter distribution is 
strongly negatively associated (CNA+). In particular, this 
solves an open problem of Pemantle \cite{P} and Liggett \cite{L2} stating 
that the
distribution of a symmetric exclusion process at time $t\geq 0$ with 
non-random/deterministic initial configuration is negatively associated, and 
shows that 
the same is actually true whenever the initial distribution is strongly 
Rayleigh. In a later paper \cite{L3}, Liggett has applied these results to 
prove convergence to the normal and Poisson laws for various functionals
of the symmetric exclusion process. 

In \S \ref{Almost} we establish
equivalences between several negative dependence properties for the class of
almost symmetric/exchangeable measures that we define in \S \ref{s21}. This 
extends Pemantle's corresponding theorem for symmetric measures and 
confirms his conjecture on strong negative
association (CNA+) in the almost exchangeable case. The results in 
\S \ref{Almost} 
are also useful in \S \ref{s-countex}, as they allow us to show that the 
examples we construct
there provide counterexamples to the series of conjectures on ultra 
log-concave rank sequences proposed in \cite{P,W} that we already alluded to. 
We note that Corollary \ref{RandC} in \S \ref{Almost} was subsequently 
proved by different methods in \cite{KN}, where it was additionally shown that 
almost exchangeable measures satisfy the (strong) Feder-Mihail property.

The aforementioned connections between
the negative dependence properties satisfied by strongly Rayleigh measures
and the geometry of zero sets of their generating polynomials are much in
the same spirit as e.g.~Rota's philosophy about the ubiquity of zeros of 
polynomials in combinatorics \cite{rota}. In the present 
context, these connections actually prove to have fruitful consequences for 
both the theory of negative dependence and the theory of 
real stable and hyperbolic polynomials, as they give new insight and further 
potential applications in these and related areas.

\section{A Plethora of Negative Dependence Properties 
and Conjectures}\label{s2-plet}

In this section we explain in terms of generating 
polynomials the probabilistic notions and negative 
dependence properties studied in e.g.~\cite{L2,P,W}, and we clarify the 
connections between these as well as 
new classes of probability measures that we introduce 
below. We then formulate the problems and conjectures 
made in \cite{L2,P,W} that we solve in the next sections.

\subsection{Negative Dependence Concepts}\label{s21}

Recall that a (real or complex) multivariate polynomial is said to be 
{\em multi-affine} 
if it has degree at most one in each variable. 
Denote by $\fP_n$, $n\in\NN$, the set of all 
probability measures on the Boolean algebra
$2^{[n]}$. (In the finite case we will always assume that the 
$\sigma$-algebra is the full algebra.) For $i\in [n]$ the $i$-th coordinate 
function on $2^{[n]}$ 
is an atomic (binary) random variable given by $X_i(S)=1$ if 
$i\in S$ and $0$ otherwise, where $S\subseteq [n]$, while the 
characteristic function $\chi_S$ of $S$ is defined by $\chi_S(T)=1$ if 
$T= S$ and $0$ otherwise.
Using the inclusion-exclusion principle one can show that any such 
characteristic function may be written as a multi-affine polynomial in 
$X_1,\ldots,X_n$. Any scalar function 
$F=\sum_{S\subseteq [n]}F(S)\chi_S$ on $2^{[n]}$ may therefore 
be viewed as a multi-affine polynomial in 
$X_1,\ldots,X_n$. By abuse of notation the latter is written as 
$F(X_1,\ldots,X_n)$ and then $\mu(F)=\int Fd\mu$. However, in what follows
the $X_i$'s always stand for (probabilistic) random variables and we will 
use other variables (such as $x,y,z,w$) when working with arbitrary 
functions/polynomials.

Let $\PP_n$ be the set of all multi-affine polynomials
in $n$ variables $f(z_1,\ldots,z_n)$ with non-negative coefficients such that
$f(\ben)=1$, where $\ben=(1,\ldots,1)\in\RR^n$ denotes the 
``all ones vector''. There is a 1-1 
correspondence between $\fP_n$ and $\PP_n$:
if $\mu\in\fP_n$ we may form its generating polynomial, namely 
$$
g_\mu(z) = \int z^{S} d\mu(S) = \sum_{S \subseteq [n]} \mu(S) z^S, \quad
z=(z_1,\ldots, z_n), \quad z^S:= \prod_{i \in S}z_i,
$$
and if 
\begin{equation}\label{coeff}
f(z)=\sum_{S\subseteq [n]}a_Sz^S\in\PP_n
\end{equation}
we define a measure $\mu_f$ on $2^{[n]}$ by setting $\mu_f(S)=a_S$, 
$S\subseteq [n]$. It is clear that $g_\nu\in\PP_n$, $\mu_f\in\fP_n$, 
$g_{\mu_f}=f$, and $\mu_{g_\nu}=\nu$ for any $\nu\in\fP_n$, $f\in\PP_n$. 
For convenience, we will sometimes use the symbol $\pa_i$ for $\pa/\pa z_i$
$1\le i\le n$, and the multi-index notation 
$\pa^S=\prod_{i\in S}\pa_i$, $S\subseteq [n]$. 

The NLC property for measures in $\fP_n$ 
translates into the corresponding property -- again called the 
{\em negative lattice condition} and denoted by NLC -- for polynomials in 
$\PP_n$: a polynomial $f\in\PP_n$ satisfies NLC if and only if its 
corresponding measure $\mu_f\in\fP_n$ does. This amounts to saying that
\begin{equation*}
\pa^S f(0,\ldots,0)\,\pa^T f(0,\ldots,0)\ge \pa^{S\cup T} 
f(0,\ldots,0)\,\pa^{S\cap T} f(0,\ldots,0), \tag*{(NLC)}
\end{equation*}
which in notation \eqref{coeff} becomes  
$a_Sa_T\ge a_{S\cup T}a_{S\cap T}$ for all $S,T\subseteq [n]$.

The so-called {\em ``closure'' operations}
discussed in e.g.~\cite[II.A]{P} may also be reformulated in terms of 
generating polynomials:

\begin{itemize}
\item[(i)] 
{\em Products.} If $\mu_1\in\fP_{n_1}$ has generating
polynomial $g_{\mu_1}(z_1,\ldots,z_{n_1})\in\PP_{n_1}$ and 
$\mu_2\in\fP_{n_2}$ has generating
polynomial $g_{\mu_2}(z_1,\ldots,z_{n_2})\in\PP_{n_2}$ then the product
$\mu_1\times \mu_2$ is the measure $\mu\in\fP_{[n_1+n_2]}$ with generating
polynomial
$g_\mu(z_1,\ldots,z_{n_1+n_2})
=g_{\mu_1}(z_1,\ldots,z_{n_1})g_{\mu_2}(z_{n_1+1},\ldots,z_{n_1+n_2})
\in\PP_{n_1+n_2}$.
\item[(ii)]
{\em Projections.} Given $S\subseteq [n]$ and $\mu\in\fP_n$ the 
projection of
$\mu$ onto $2^S$ is the measure $\mu'\in\fP_{|S|}$ with generating polynomial
$$g_\mu(z_1,\ldots,z_n)\big|_{z_i=1,\,i\in [n]\setminus S}\in\PP_{|S|}.$$
Note that in order for this definition 
to make sense we should also relabel the variables using indices in $[|S|]$. 
However, we allow ourselves this abuse of notation here and in what follows. 
\item[(iii)]
{\em Conditioning.} Let $\mu\in\fP_n$ with generating polynomial
$g_\mu(z_1,\ldots,z_n)\in\PP_n$ and fix some $i\in [n]$. The measure 
obtained from $\mu$ by conditioning the $i$-th random variable $X_i$ 
to be $0$ is the measure on $2^{[n]\setminus\{i\}}$ with generating polynomial
$$\frac{g_\mu(z_1,\ldots,z_n)\big|_{z_i=0}}
{g_\mu(z_1,\ldots,z_n)\big|_{z_i=0,\,z_j=1,\,j\neq i}}\in\PP_{n-1}$$
while the measure obtained from $\mu$ by conditioning $X_i$ 
to be $1$ is the measure on $2^{[n]\setminus\{i\}}$ with generating polynomial
$$\lim_{z_i\to\infty}\frac{g_\mu(z_1,\ldots,z_n)}
{g_\mu(z_1,\ldots,z_n)\big|_{z_j=1,\,j\neq i}}
=\frac{\pa_i g_\mu(z_1,\ldots,z_n)}
{\pa_i g_\mu(1,\ldots,1)}\in\PP_{n-1}.$$
Note that these conditional probability measures are well defined provided 
that the denominators appearing in the above expressions are non-zero.
\item[(iv)]
{\em External fields} (as pointed out in \cite[II.A]{P}, this name is 
borrowed 
from the Ising model). If  $\mu\in\fP_n$ has generating polynomial
$g_\mu(z_1,\ldots,z_n)\in\PP_n$ and $a_i$, $1\le i\le n$, are non-negative real
numbers then the measure obtained from $\mu$ by imposing the external field
$(a_1,\ldots,a_n)$ is the measure in $\fP_n$ with generating polynomial
$$\frac{g_\mu(a_1z_1,\ldots,a_nz_n)}{g_\mu(a_1,\ldots,a_n)}\in\PP_n,$$
which is well defined provided that $g_\mu(a_1,\ldots,a_n)\neq 0$.
\item[(v)]
{\em Symmetrization.} Denote the symmetric group on $n$ 
elements by $\mathfrak{S}_n$. For $\sigma\in \mathfrak{S}_n$ and 
$S\subseteq [n]$ let
$\sigma(S)=\{\sigma(s): s \in S\}$. Given $\mu\in\fP_n$ define a measure
$\sigma(\mu)\in\fP_n$ by setting $\sigma(\mu)(S)=\mu(\sigma(S))$, 
$S\subseteq [n]$. The (full or complete) symmetrization of 
$\mu$ is the measure 
$$\mu_s=\frac{1}{n!}\sum_{\sigma\in\mathfrak{S}_n}\sigma(\mu)\in \fP_n.$$ 
Note that the generating polynomial of $\mu_s$ is
\begin{equation}\label{sym}
\begin{split}
g_{\mu_s}(z)=\frac{1}{n!}\sum_{\sigma\in\mathfrak{S}_n}
g_{\sigma(\mu)}(z)
&=\frac{1}{n!}\sum_{\sigma\in\mathfrak{S}_n}\sum_{S\subseteq [n]}
a_{\sigma(S)}z^S\\
&=\sum_{k=0}^{n}\frac{\sum_{|S|=k}a_S}{\binom{n}{k}}
e_k(z_1,\ldots,z_n)\in\PP_n,
\end{split}
\end{equation}
where $\sum_{S\subseteq [n]}a_{S}z^S=g_{\mu}(z)$ and $e_k$, $0\le k\le n$, 
is the $k$-th elementary symmetric function (see \cite[\S 3.4]{W} for a proof
of the last identity). 
\item[(vi)]
{\em Partial symmetrization.} Let $\mu\in\fP_n$, $0\leq \theta \leq 1$, 
$1\leq i<j \leq n$, and $\tau=(ij)\in\mathfrak{S}_n$ be the transposition that 
permutes $i$ and $j$. The partial symmetrization of $\mu$ with respect to
$\tau$ and $\theta$ is the measure 
\begin{equation}\label{t-theta}
\mu^{\tau,\theta}
=\theta\mu+(1-\theta)\tau(\mu)\in\fP_n
\end{equation}
whose generating polynomial is obviously 
$g_{\mu^{\tau,\theta}}=\theta g_{\mu}+(1-\theta)g_{\tau(\mu)}\in\PP_n$.
\end{itemize}

\begin{remark}\label{r-complex}
All these ``closure'' operations are in fact well defined for arbitrary 
complex measures on $2^{[n]}$ if one drops the normalization factors 
appearing in the denominators of the
above expressions (these were used just to make sure that the resulting 
measures are again probability measures). 
In \S \ref{s-stable} we will actually use the partial symmetrization 
procedure \eqref{t-theta} for complex measures on $2^{[n]}$.
\end{remark}

The corresponding ``closure'' operations for generating polynomials -- i.e., 
{\em products}, 
{\em projections}, {\em conditioning}, {\em external fields}, 
{\em symmetrization}, {\em partial symmetrization} and {\em truncations} 
(cf.~Definition \ref{d-trunc} in \S \ref{s24})-- 
are defined simply by considering the resulting 
polynomials for each of the operations discussed above for measures.

\begin{definition}\label{d-exc}
A complex measure $\mu$ on $2^{[n]}$ is called {\em symmetric} or
{\em exchangeable} if $\sigma(\mu)=\mu$ for any $\sigma\in \mathfrak{S}_n$,
or equivalently, its generating polynomial $g_\mu$ is symmetric in all $n$
variables. We say that $\mu$ is {\em almost symmetric} or {\em almost
exchangeable} if $g_\mu$ is symmetric in all but possibly one variable. 
\end{definition}

\begin{definition}\label{d-inc-ev}
If $\mathcal{A}$ is a collection of subsets of $[n]$ we let 
$\chi_{\mathcal{A}} : 2^{[n]} \rightarrow \RR$ be the characteristic function 
of $\mathcal{A}$ defined by 
$\chi_{\mathcal{A}}(S) = 1$ if $S \in \mathcal{A}$ and 
$\chi_{\mathcal{A}}(S) = 0$ if $S \notin \mathcal{A}$.
An {\em increasing event} $\mathcal{A}$ on $2^{[n]}$ is a collection of 
subsets of $[n]$ that is closed upwards under containment, i.e., if 
$A \in \mathcal{A}$ and $A \subseteq B \subseteq [n]$ then 
$B \in \mathcal{A}$. Such an event {\em depends} only on the set
$\cI(\cA):=I_1\cup\cdots\cup I_m$ (that is, on the variables/coordinate 
functions $X_i$ with 
$i\in\cI(\cA)$), where $I_1,\ldots,I_m$ are the minimal sets of $\cA$ 
with respect to inclusion. Any (non-identically zero) non-negative increasing
function $f$ on $2^{[n]}$ may be written as $f=\sum_{i=1}^{k}a_i\chi_{\cA_i}$ 
for some increasing events $\cA_i$ on $2^{[n]}$ and $a_i>0$, $1\le i\le k$.
Clearly, $f$ depends on the set $\cI(f):=\cI(\cA_1)\cup\cdots\cup\cI(\cA_k)$.
\end{definition}

The weakest negative dependence property is the following 
(see, e.g., \cite{P}).

\begin{definition}\label{d-pnc}
A measure $\in\fP_n$ is said to be {\em pairwise
negatively correlated} or p-NC for short 
if $\mu(X_i)\mu(X_j)\ge \mu(X_iX_j)$ whenever 
$1\le i\neq j\le n$. In terms of the generating polynomial $g_\mu$ this 
translates into the inequalities
$$\pa_i g_\mu(\ben)\pa_j g_\mu(\ben)\ge \pa_i \pa_j g_\mu(\ben),
\quad 1\le i\neq j\le n,$$
where $\ben\in\RR^n$ is as before the ``all ones vector''. We say
that a polynomial $f\in\PP_n$ is p-NC (pairwise negatively correlated) if it
satisfies the above inequalities. 
\end{definition}

As shown by the following example, NLC does not necessarily imply p-NC.

\begin{example}\label{ex1}
The measure $\mu\in\fP_4$ with generating polynomial
$$g_\mu(z_1,z_2,z_3,z_4)=\frac{1}{2}(z_1z_2+z_3z_4)$$
trivially satisfies NLC but $\mu(X_1X_2)=\frac{1}{2}$ while 
$\mu(X_1)\mu(X_2)=\frac{1}{4}$, so $\mu$ is not p-NC. 
\end{example}

\begin{definition}[Definition 2.4 in \cite{P}]\label{h-nlc}
A measure $\mu\in\fP_n$ or a polynomial $f\in\PP_n$ 
satisfies the {\em hereditary negative lattice condition} or $\hNLC$ 
if every projection satisfies $\NLC$. One further says that $\mu$ 
(respectively, 
$f$) satisfies the {\em strong  hereditary negative lattice 
condition} or $\hNLC+$ if every measure (respectively, polynomial) obtained 
from $\mu$ (respectively, $f$) by imposing an external field satisfies
h-NLC.
\end{definition}

\begin{definition}\label{d1}
A polynomial $f\in\PP_n$ is called a {\em Rayleigh polynomial} if
\begin{equation}\label{ral-ineq}
\frac {\partial f} {\partial z_i}(x) \frac {\partial f} {\partial z_j}(x)
\geq \frac {\partial^2 f} {\partial z_i \partial z_j}(x) f(x)
\end{equation}
for all $x=(x_1,\ldots,x_n)\in \RR_+^n$ and $1\leq i,j \leq n$, where as 
usual $\RR_+=(0,\infty)$. More 
generally, a multi-affine 
polynomial in $\RR[z_1, \ldots, z_n]$ with non-negative coefficients is 
called a {\em Rayleigh polynomial} if it satisfies the above condition. 
A measure $\mu\in\fP_n$ is said to be a {\em Rayleigh measure} if its
generating polynomial $g_\mu$ is Rayleigh.
\end{definition}

The notion of Rayleigh polynomial was introduced in \cite[\S 3.1]{W}, the 
terminology being motivated by its similarity with the Rayleigh monotonicity 
property of 
the (Kirchhoff) effective conductance of linear resistive electrical networks,  
see \cite{CW,W2}. It was first considered for uniform measures on the set of bases of a matroid \cite{CW} as a strengthening of weaker notions studied in \cite{feder,seymour-welsh}. 
The following properties of
Rayleigh polynomials are consequences of Definition~\ref{d1}.
\begin{proposition}\label{p-proper}
If $f(z_1,\ldots,z_n)$ is Rayleigh then so are the following polynomials:
\begin{enumerate}
\item $\pa^S f(z_1,\ldots,z_n)$ for any $S\subseteq [n]$;
\item $f(z_1+\al_1,\ldots,z_n+\al_n)$ whenever $\al_i\ge 0$, $i\in [n]$;
\item $f(z_1,\ldots,z_n)|_{z_i=\al_i}$ for any $\al_i\ge 0$, $1\le i\le n$; 
\item $f(a_1z_1,\ldots,a_nz_n)$ for all $a_i\ge 0$, $i\in [n]$;
\item the ``inversion'' of $f$, i.e., the polynomial 
$z_1\cdots z_nf(z_1^{-1},\ldots,z_n^{-1})$;
\item the polynomial in $nk$ variables 
$f\!\left(k^{-1}\sum_{i=1}^k z_{1i},\ldots,k^{-1}\sum_{i=1}^k z_{ni}\right)$.
\end{enumerate}
\end{proposition}

To prove (1) note that by induction it is enough to do it for $S=\{i\}$ and
arbitrary $i\in [n]$. This follows by first writing 
$$f(z_1,\ldots,z_n)=z_i\pa_i f(z_1,\ldots,z_n)+f(z_1,\ldots,z_n)|_{z_i=0}$$
and then letting $z_i\to\infty$ in the Rayleigh inequalities \eqref{ral-ineq}
for $f$ corresponding to pairs of 
variables indexed by distinct $j,k\in [n]\setminus \{i\}$.

\begin{proposition}\label{p-r-h}
A measure in $\fP_n$ or a polynomial in $\PP_n$ is Rayleigh
if and only if it is h-NLC+.
\end{proposition}

\begin{proof}
In \cite[Theorem~4.4]{W} it was proved that if $f\in\PP_n$ is Rayleigh then 
$f$ satisfies NLC, which combined with Proposition~\ref{p-proper} shows that
$f$ is h-NLC+. To prove the converse statement, let $i,j\in [n]$ with $i\neq j$
and set 
$$g(z_i,z_j)=f(z_1,\ldots,z_n)|_{z_k=1,\,k\in [n]\setminus \{i,j\}}.$$
If $f$ is h-NLC+ then $g$ satisfies NLC (which in this case implies that
$g(z_i,z_j)$ is in fact a real stable polynomial, see 
Definition~\ref{d-rs-pol} in \S \ref{s23} and Theorem~\ref{stable-Rayleigh}
in \S \ref{s-stable}). In particular, we deduce that 
$\pa_i \pa_j f(\ben)\le \pa_i f(\ben)\pa_j f(\ben)$, 
where $\ben$ stands for the ``all ones vector'' as in Definition~\ref{d-pnc}.
Using the external field condition given in Definition~\ref{h-nlc}
wee see that the aforementioned inequality also holds for the polynomial
$f(x_1z_1,\ldots,x_nz_n)$ whenever $x_i\ge 0$, $i\in [n]$, which by 
Definition~\ref{d1} amounts to saying that $f$ is Rayleigh.
\end{proof}

\begin{remark}\label{r-hnlc}
It is interesting to note that one actually has the following analog of 
Proposition \ref{p-r-h}
for the h-NLC property: a measure $\mu\in\fP_n$ -- or its generating polynomial
$g_\mu\in\PP_n$ -- is h-NLC if and only if $g_\mu$ satisfies the Rayleigh 
inequalities
\eqref{ral-ineq} for all $x=(x_1,\ldots,x_n)$ with $x_i\in\{0,1,\infty\}$,
$i\in [n]$.
\end{remark}

\begin{definition}\label{phr}
A {\em homogeneous Rayleigh measure} is one whose generating polynomial
is homogeneous and Rayleigh. The set of all measures that are projections of 
homogeneous Rayleigh measures is denoted by PHR.
\end{definition}

Let us now recall the key concept of negative association. For probability 
measures this notion is usually defined as follows (cf., e.g., \cite{L2,P,W}).

\begin{definition}\label{d-na}
A measure $\mu\in\fP_n$ is called {\em negatively associated} or $\NA$ 
if
$$\int F d\mu \int G d\mu\ge \int FG d\mu$$
for any increasing functions $F,G$ on $2^\cN$ that depend on 
disjoint sets of coordinates (cf.~Definition~\ref{d-inc-ev}). One says that 
$\mu$ is {\em conditionally negatively associated}
or $\CNA$ if each measure obtained from $\mu$ by conditioning on some (or none)
of the values of the variables is $\NA$. Finally, $\mu$ is called {\em strongly
conditionally negatively associated} or $\CNA+$ if each measure obtained from 
$\mu$ by imposing external fields and projections is $\CNA$.
\end{definition}

\begin{remark}\label{r-pair}
It is clear from the definitions
that each of the five properties h-NLC, Rayleigh/h-NLC+, NA, CNA, CNA+ implies 
p-NC.
\end{remark}

\begin{remark}\label{r-cn-hn}
As explained in \cite{P}, one has the following subordination relations: 
CNA $\Rightarrow$ h-NLC and CNA+ $\Rightarrow$ Rayleigh/h-NLC+. In 
\S \ref{s-sna} and \S \ref{s-countex} below we show that PHR 
$\Rightarrow$ CNA+ and that this implication is strict.
\end{remark}

We note that NA can fail to imply NLC even for symmetric measures:

\begin{example}\label{ex2}
Let as before $e_k$, $0\le k\le n$, be the $k$-th elementary symmetric function
in $n$ variables and consider the measure $\mu\in\fP_3$ with 
generating polynomial
$$
g_\mu(z_1,z_2,z_3)=\frac{1}{15}\left[3+2e_1(z_1,z_2,z_3)+2e_2(z_1,z_2,z_3)
\right].
$$
It is not difficult to check that $\mu$ is NA but not NLC.
\end{example}

\begin{definition}\label{d-ulc}
The {\em diagonal 
specialization} of a polynomial $f\in\CC[z_1,\ldots,z_n]$ is the univariate 
polynomial 
$$t\mapsto\De(f)(t):=f(t,\ldots,t).$$
A real sequence $\{a_k\}_{k=0}^n$ is {\em log-concave} 
if $a_k^2 \geq a_{k-1}a_{k+1}$, $1\leq k \leq n-1$, and it is said to have
{\em no internal zeros} if the indices of its non-zero terms form an interval
(of non-negative integers). 
We say that a 
{\em non-negative} sequence $\{a_k\}_{k=0}^n$ is
\begin{itemize}
\item $\LC$ if it is log-concave with no internal zeros; 
\item $\SLC$ (strongly log-concave) if the sequence 
$\{k!a_k\}_{k=0}^{n}$ is $\LC$;
\item $\ULC$ (ultra log-concave) if the sequence
$\left\{a_k/\binom{n}{k}\right\}_{k=0}^{n}$ is $\LC$. 
\end{itemize}
Clearly, $\ULC\Rightarrow \SLC\Rightarrow 
\LC$. If $\mu\in\fP_n$ the sequence 
$$\left\{\mu\left(\sum_{i=1}^{n}X_i=k\right)\right\}_{k=0}^{n}
=\left\{\frac{\De(g_\mu)^{(k)}(0)}{k!}\right\}_{k=0}^{n}$$
is called the {\em rank sequence} of $\mu$ as well as of its generating 
polynomial $g_\mu\in\PP_n$. A measure in $\fP_n$ is then said to be $\ULC$,
$\SLC$ or $\LC$ if its rank sequence is $\ULC$,
$\SLC$ or $\LC$, respectively. 
\end{definition}

\subsection{Strongly Rayleigh Measures}\label{s23}
Motivated by similar notions in the theory of multivariate entire functions
studied by Levin \cite{Le}, the following definition was made in
\cite{BBS,BBS1,BBS2} (see also \cite{BB,BB-2}).

\begin{definition}\label{d-rs-pol}
A polynomial $f \in \CC[z_1,\ldots, z_n]$ is called {\em stable} if 
$f(z_1,\ldots, z_n) \neq 0$ 
whenever $\Im(z_j) >0$ for $1\le j\le n$. A stable polynomial with all
real coefficients is called {\em real stable}.
\end{definition}

Multivariate stable polynomials are intimately connected with other fundamental
objects such as hyperbolic polynomials, Lee-Yang polynomials and 
polynomials with the half-plane property. Polynomials of the aforementioned 
types are widely encountered in both mathematics and physics
\cite{ABG,BBCV,BBS,BBS1,BBS2,BB,PB,COSW,guler,garding,LS} (see also 
\S \ref{s-gar}). 

Let us define a class of measures induced by stable polynomials.

\begin{definition}\label{def-rsm}
A measure $\mu\in\fP_n$ is called {\em strongly Rayleigh} if its generating 
polynomial $g_\mu\in\PP_n$ is (real) stable. For convenience, we will  
sometimes refer to a real stable polynomial $f\in\PP_n$ as a 
{\em strongly Rayleigh polynomial}.
\end{definition}

\begin{remark}\label{r-term}
As we explain in \S \ref{ss-real-s} and \S \ref{s-stable}, a strongly 
Rayleigh measure is necessarily
Rayleigh. The terminology adopted is 
motivated by the 
fact that the equivalent condition for real stability for multi-affine 
polynomials in
Theorem~\ref{stable-Rayleigh} is a natural strengthening of the definition 
of a 
Rayleigh polynomial/measure (Definition~\ref{d1}). The term is taken 
from \cite{CW}, where a 
matroid was defined to be strongly Rayleigh if \eqref{stable} holds for
 its bases-generating polynomial. However, it was not known 
then that condition \eqref{stable} is equivalent to stability. 
\end{remark}

A simple albeit important example of strongly Rayleigh measure is the
following.

\begin{definition}\label{d-prod-m}
A {\em product measure} on $2^{[n]}$
is a measure $\mu\in\fP_n$ with generating polynomial $g_\mu$ of the form
$$g_\mu(z_1,\ldots,z_n)
=\prod_{i=1}^{n}(q_iz_i+1-q_i),\quad 0\leq q_i\leq 1,\,i\in [n].$$
One says that $\mu$ is
{\em deterministic} or {\em non-random} if it is a point mass
on $2^{[n]}$, i.e., if there exists $S\subseteq [n]$ such that
$q_i=1$ if $i\in S$ and $q_j=0$ for  
$j\in [n]\setminus S$.
\end{definition}

In \S \ref{s3-ex} we show that the class of strongly Rayleigh measures 
contains numerous other 
important examples of measures appearing in matrix theory,
combinatorics, probability theory, particle systems, 
and in \S \ref{s-stable} we study this class in detail. 

\subsection{Symmetric Homogenization}\label{s22}
We will now define a symmetric homogenization 
procedure for arbitrary measures in $\fP_n$. This natural construction 
proves to be quite useful for studying strongly Rayleigh measures 
(\S \ref{s-gar}).

\begin{definition}\label{d-sh}
Given a measure $\mu\in\fP_n$ define a new measure $\mu_{sh}\in\fP_{2n}$ 
called the {\em symmetric homogenization} of $\mu$ by
$$
\mu_{sh}(S) = 
\begin{cases}
\mu(S\cap[n])\binom{n}{|S\cap [n]|}^{-1} \mbox{ if } |S|=n, \\
0 \mbox{ otherwise. } 
\end{cases} 
$$
\end{definition}

Note that the generating polynomial of $\mu_{sh}$ may be written as
$$g_{\mu_{sh}}(z_1,\ldots,z_n,z_{n+1},\ldots,z_{2n})=
\sum_{S \subseteq [n]} \mu(S)\binom{n}{|S|}^{-1} z^S 
e_{n-|S|}(z_{n+1},\ldots, z_{2n}),
$$
where $z=(z_1,\ldots,z_n)$ and $e_k$, $0\le k\le n$, denotes the $k$-th 
elementary symmetric function. Clearly, $g_{\mu_{sh}}$ is symmetric in the 
variables $z_{n+1},\ldots,z_{2n}$ and $\mu_{sh}$ projects to $\mu$, that is,
$g_{\mu_{sh}}(z_1,\ldots,z_n,1,\ldots,1)=g_{\mu}(z_1,\ldots,z_n)$. It is not 
difficult to show that 
$g_{\mu_{sh}}$ is actually the unique polynomial
in $\PP_{2n}$ of total degree $n$ that is 
symmetric in $z_{n+1},\ldots,z_{2n}$ and satisfies the latter identity. 
Therefore, $\mu_{sh}$ is canonically 
determined by the aforementioned conditions.

\subsection{Measures on countably infinite sets}\label{s240}
In \S \ref{n-s6} we will need to extend some of the notions of  
negative dependence to probability 
measures on $2^{E}$, where $E$ is a countably infinite set. This is 
canonically done as follows. 

\begin{definition}\label{infinite}
Let $\mu$ be a probability measure on $2^E$, where $E$ is a countably 
infinite set satisfying 
$$
\mbox{ (A): \ \ }\mbox{ for all } x \in E \mbox{ the set } 
\left\{T \in 2^E : x \in T\right\} \mbox{ is measurable. } 
$$
 Let $\mathcal{P}$ be a property defined for measures on power-sets of finite 
sets. We say that $\mu$ has property $\mathcal{P}$ provided that the 
projection of $\mu$ to all finite subsets of $E$ has property $\mathcal{P}$. 
(Note that 
condition (A) is necessary for projections on finite sets to be well-defined.) 
\end{definition}

\subsection{Towards a Theory of Negative Dependence: Problems and 
Solutions}\label{s24}

A series of problems and conjectures was proposed by several authors as 
natural steps toward building as general a theory of negative dependence as 
possible. In \cite[\S III]{P} Pemantle pointed out that although
they are desirable, the many negative dependence properties listed in 
\S \ref{s21} might not be mutually satisfiable, and he suggested a problem that
may be formulated as follows.

\begin{problem}\label{pem-prob}
Find a natural and useful class of ``negatively dependent measures''.
\end{problem}

Note that there may well exist several such classes, if indeed any, 
depending on the
applications that one has in mind. We show that the class of strongly Rayleigh 
measures introduced in \S \ref{s23} enjoys all the negative dependence 
properties defined in \S \ref{s21} as well as other properties, which 
provides an answer to Problem~\ref{pem-prob}. 

The problem of describing natural
negative dependence properties that are preserved by symmetric exclusion 
evolutions has 
attracted some attention in the theory of 
interacting particle systems
and Markov processes \cite{L0,L00,L2}, see Problem~\ref{p-see} in 
\S \ref{ss-ex-case} and the discussion therein. In particular, the 
following conjecture
was implicitly made in \cite[p.~1374]{P} and \cite[p.~551]{L2}. 

\begin{conjecture}\label{con-lig}
If the initial configuration of a symmetric exclusion process
is non-random/deterministic, then the distribution at time $t$ is NA for 
all $t>0$.
\end{conjecture}

In \S \ref{n-s6} we prove a stronger form of this conjecture using the 
results on strongly Rayleigh measures obtained in \S \ref{s-stable}.

There are several conjectures, first published in \cite{P} 
but apparently 
of obscure folklore origin (cf.~\cite{W}), stating that various negative
dependent properties give rise to ULC rank sequences. The strongest form
of these is contained in the first part of Conjecture 4 in Pemantle's paper 
\cite{P} (see also \cite[Conjecture 3.9]{W}) and reads as follows.

\begin{conjecture}\label{c-p}
If $\mu\in\fP_n$ is $\NA$ then it is $\ULC$.
\end{conjecture}

Inspired by Conjecture~\ref{c-p} 
and motivated by Mason's long-standing
conjecture in enumerative graph/matroid theory \cite{Mason}, Wagner made the 
following conjecture in \cite[Conjecture 3.4]{W} and baptized it 
``The Big Conjecture'' since -- if true -- it would prove the validity of 
Mason's conjecture for a large class of matroids.  

\begin{conjecture}\label{big}
Any Rayleigh measure $\mu\in\fP_n$ is $\ULC$.
\end{conjecture}

Since we now know that the Rayleigh and h-NLC+ properties are equivalent 
(Proposition~\ref{p-r-h}), 
Conjecture~\ref{big} is actually already contained in the second part
of Conjecture 4 in Pemantle's paper \cite{P}. The weakest form of the latter 
conjecture is the one stated below.

\begin{conjecture}\label{cna-ulc}
If $\mu\in\fP_n$ is $\CNA+$ then it is $\ULC$.
\end{conjecture}

In \S \ref{s-countex} we construct examples of measures in $\fP_n$ for any
$n\ge 20$ that disprove all these three conjectures (for the strongest
of these, namely Conjecture~\ref{c-p}, we can actually find
counterexamples already in $\fP_3$ and $\fP_4$, see Remark~\ref{r-strong-NA}).
In fact, our examples 
show that the strongest negative association property (CNA+) 
need not even give rise to an SLC rank sequence 
(Definition~\ref{d-ulc}), thus invalidating a still weaker version of 
Conjecture~\ref{cna-ulc}. Moreover, these same examples show that
neither the PHR property (Definition~\ref{phr}) nor
indeed any of the properties weaker than strongly Rayleigh implies the 
ULC property, and also that Conjecture~\ref{cna-ulc} fails even under extra 
assumptions like those suggested in \cite[Problem 6]{BBCV}. This confirms 
once again the delicate nature of negative dependence.

The FKG theorem \cite{FKG} is a powerful tool that allows one to establish 
(global) positive correlation inequalities and limit theorems from a 
local (usually easier to check) condition, namely the positive lattice 
condition PLC. By contrast, there is as yet no general and practical 
local-to-global device for negative dependence/association. 
The NLC property is not closed under 
projections and implies neither the NA nor the (weaker) p-NC property, so the
negative version of the FKG theorem fails. As a potential remedy to this
situation, in \cite[Conjecture~2]{P} Pemantle
conjectured that the implications 
CNA $\Rightarrow$ h-NLC and CNA+ $\Rightarrow$ Rayleigh/h-NLC+ 
(established in \cite{P}, cf.~Remark~\ref{r-cn-hn}) are actually equivalences:

\begin{conjecture}\label{con-cna-R}
h-NLC $\,\Rightarrow$ CNA and Rayleigh/h-NLC+ $\Rightarrow$ CNA+. 
\end{conjecture}

We prove Conjecture~\ref{con-cna-R} 
for almost symmetric measures (cf.~Definition~\ref{d-exc}). 
Since PHR measures are CNA+  
(Theorem~\ref{PHR-NA}), we are naturally
led to consider the following related problem. 

\begin{problem}\label{prob-phr-cna}
Is PHR $=$ Rayleigh/h-NLC+?
\end{problem}

Note that if true, this would prove the second part of 
Conjecture~\ref{con-cna-R}. 
In \S \ref{s-countex} we show that the answer to 
Problem~\ref{prob-phr-cna} is negative, and that in fact the class of 
PHR measures 
is strictly contained in the class of CNA+ measures. The question 
whether Conjecture~\ref{con-cna-R} is true in full generality 
therefore remains open.

We summarize the implications between the various notions of 
negative dependence discussed so far in Fig.~\ref{fig} below. 
Compared with
the one in \cite[\S II.B]{P}, this diagram contains five extra 
classes -- three ``old'' (ULC, p-NC, NLC) and two ``new'' 
(strongly Rayleigh, PHR) -- and disregards the classes called 
JNRD, JNRD+ in \cite{P}.

\begin{figure}[!htb]
\xymatrix{
& & & \text{ULC} & {\text{{\normalsize Strongly}}\atop 
\text{{\normalsize Rayleigh}}}\ar@{>}[l]\ar@{>}[d] & \\
& & & & \text{PHR}\ar@{>}[d] & \\
& & & & \text{CNA+}\ar@{>}[d]\ar@{>}[r] & {\text{{\normalsize Rayleigh}}\atop 
\text{{\normalsize (h-NLC+)}}}\ar@{>}[d] \\
& & & & \text{CNA}\ar@{>}[d]\ar@{>}[r] & \text{h-NLC}\ar@{>}[ddl]\ar@{>}[d] \\
& & & & \text{NA}\ar@{>}[d] & \text{NLC} \\
& & & & \text{p-NC} & }
\caption{Subordination relations between negative dependence properties.}
\label{fig}
\end{figure}
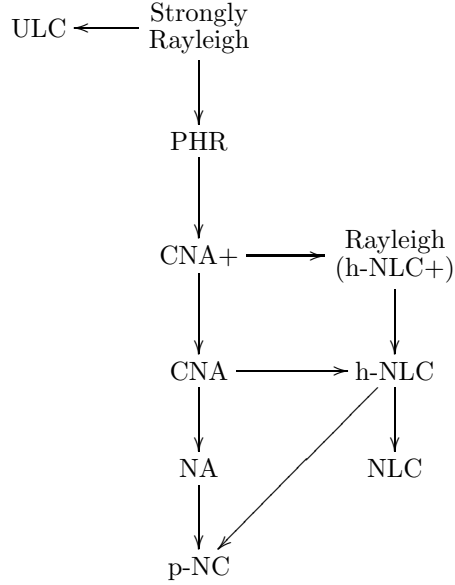

As we already noted, the above examples and those
in the next sections show that the 
left horizontal, down left, and all vertical implications 
in Fig.~\ref{fig} are strict, while the question whether the right horizontal
implications are also strict is open (cf.~Conjecture~\ref{con-cna-R}). Note
also that there is no arrow pointing towards ULC from anything ``below'' 
strongly Rayleigh.

In \cite[\S I.B \& \S III.C]{P} it was pointed out that the notions of 
stochastic domination and truncation are also 
useful when studying (positive or) negative dependence.

\begin{definition}\label{d-stoc}
If $\mu,\nu\in\fP_n$ are such that $\mu(\cA)\ge \nu(\cA)$ for any
increasing event $\cA$ on $2^{[n]}$ one says that $\mu$ 
{\em stochastically dominates} $\nu$, written $\mu\suc \nu$ or $\nu\pc \mu$.
\end{definition}

\begin{definition}\label{d-trunc}
Let $\mu\in\fP_n$ and $1\le p\le q\le n$. The 
{\em truncation}
of $\mu$ to $[p,q]$ is the conditional measure 
$$\mu_{p,q}:=\left(\mu\,\bigg|\,p\le \sum_{i=1}^{n}X_i\le q\right)
\in\fP_n,$$
which is well defined provided that 
$\mu\left(\left\{S\in 2^{[n]}:p\le |S|\le q\right\}\right)\neq 0$. 
The generating polynomial of $\mu_{p,q}$ is thus given by
$$g_{\mu_{p,q}}(z_1,\ldots,z_n)
=\frac{g_{p,q}(z_1,\ldots,z_n)}
{g_{p,q}(1,\ldots,1)}\in\PP_n,$$
where $g_{p,q}(z)=\sum_{p\le |S|\le q}a_Sz^S$ and 
$\sum_{S\subseteq [n]}a_{S}z^S=g_{\mu}(z)$.
\end{definition}

For simplicity, we let $\mu_k=\mu_{k,k}$, $0\le k\le n$.
In \cite{J-DP,P} and references therein it has been discussed whether 
truncations of measures are operations that preserve negative dependence
properties. Question 10 of \cite{P} asks the following.

\begin{problem}\label{p-stoc-d}
Under what hypotheses on $\mu$ can one prove that $\mu_k\pc \mu_{k+1}$?
\end{problem}

Moreover, the next conjecture was made in \cite[Conjecture 8]{P}.

\begin{conjecture}\label{st}
Suppose that  $\mu$ is $\CNA+$ and that  
$$\mu(\{S \in 2^{[n]} : |S|=k\})\mu(\{S \in 2^{[n]} : |S|=k+1\})>0.$$ 
Then $\mu_k \pc \mu_{k+1}$. 
\end{conjecture}

To underline the relevance of these questions, in \cite[\S III.C]{P}
Pemantle argued that a combination of positive answers to
Problem~\ref{p-stoc-d} and Conjecture~\ref{st} (for Rayleigh 
measures) would imply Conjecture~\ref{con-cna-R}. In
\S \ref{ss-trunc} we give an answer to Problem~\ref{p-stoc-d} and prove that
its conclusion (hence also that of Conjecture~\ref{st}) 
holds for strongly Rayleigh measures, while in \S \ref{s-countex} we disprove
Conjecture~\ref{st} in its full generality.

\section{Negative Dependence in Mathematics and Physics}\label{s3-ex}

\subsection{Real Stable Polynomials and Transcendental Entire 
Functions}\label{ss-real-s}

In \cite{PB} it was proved that a multi-affine polynomial 
$f \in \RR[z_1,\ldots,z_n]$ is stable (in the sense of 
Definition~\ref{d-rs-pol}) 
if and only if 
\begin{equation}\label{stable}
\frac {\partial f} {\partial z_i}(x) \frac {\partial f} {\partial z_j}(x)
\geq \frac {\partial^2 f} {\partial z_i \partial z_j}(x) f(x)
\end{equation}
for all $x \in \RR^n$ and $1\leq i<j \leq n$ (see also 
Theorem~\ref{stable-Rayleigh} in \S \ref{s-stable} below). Therefore, 
a stable multi-affine polynomial with non-negative coefficients is 
Rayleigh/h-NLC+. As we will prove in \S \ref{s-stable}, such polynomials are 
in fact CNA+.

Clearly, by setting all variables equal in a stable polynomial one gets a 
univariate stable polynomial. Hence, if $f(z_1,\ldots,z_n)$ is stable then 
so is its diagonal specialization (see Definition~\ref{d-ulc})
$$
\De(f)(t)=\sum_{k=0}^n \binom n k c_k t^k.
$$
A univariate polynomial with real coefficients is stable if and only if it 
has all real zeros, so the familiar Newton inequalities 
(cf., e.g., \cite{CC,HLP}) imply that $f$ satisfies the conclusion of 
Conjecture~\ref{big} if $f$ is stable with non-negative coefficients. 

Multivariate generalizations of Newton's inequalities for homogeneous real
stable polynomials with non-negative coefficients were
obtained in \cite[Theorem 3]{BBS1}. Corollaries 2 and 3 of 
{\em op.~cit.~}show that the coefficients of (symmetric) such polynomials 
also satisfy unimodality properties in the sense of majorization/stochastic
domination \cite{MO}. In \S \ref{ss-trunc} we establish natural analogs
of these results for truncations of (non-homogeneous) multi-affine 
stable polynomials with non-negative coefficients.

The well-known Laguerre-Tur\'an inequalities for univariate polynomials and 
transcendental entire functions in the Laguerre-P\'olya class \cite{CC}
amount to saying that the sequence of Taylor coefficients of such a function
is SLC (Definition~\ref{d-ulc}). 
One can argue that both Newton's inequalities and the 
Laguerre-Tur\'an inequalities are 
natural manifestations of negative dependence properties encoded in 
the {\em geometry} of the zero sets of these functions. Indeed, as we will 
see in \S \ref{s-stable}, the geometric features of real stability 
(Definition~\ref{d-rs-pol}) are a key tool in our study of negative 
dependence properties for strongly Rayleigh measures and allow us to
establish many such properties for this class of measures, including:
\begin{itemize}
\item the strongest form of negative association (CNA+);
\item closure under operations (i)--(vi) defined in \S \ref{s21};
\item closure under symmetric homogenization;
\item stochastic domination properties for truncations;
\item the ULC property for rank sequences.
\end{itemize}

Rayleigh polynomials that are not stable do not seem 
to have a geometric description by means of their zero sets similar to that
of stable polynomials. 
This may explain why several negative dependence properties fail for 
arbitrary (non-strongly) Rayleigh
polynomials/measures (see \S \ref{s-countex}).

We summarize below some of the closure properties of real stable 
polynomials (see \cite{BBS1,BBS2} for their proofs) that we need
for later use.

\begin{proposition}\label{pro-rs}
Let $\HH_n(\RR)$ be the set of all real stable polynomials in $n$
variables. 
\begin{enumerate}
\item $f\in \HH_n(\RR)$ if and only if for any $\lambda\in\RR_+^n$ and 
$\mu\in\RR^n$ the univariate polynomial $t\mapsto f(\lambda t+\mu)$ has all 
real zeros;
\item  If $f\in\HH_n(\RR)$ has degree $d_j$ in variable $z_j$, $j\in [n]$,
then
\begin{itemize}
\item $\partial_j f \in \HH_n(\RR)\cup\{0\}$ for 
$1 \leq j \leq n$;
\item $f(z_1, \ldots, z_{j-1}, \alpha z_j, z_{j+1}, \ldots, z_n) 
\in \HH_{n}(\RR)$ for $j \in [n]$, $\alpha > 0$;
\item $f(z_1, \ldots, z_{j-1},\beta, z_{j+1}, \ldots, z_n) 
\in \HH_{n-1}(\RR)\cup\{0\}$ for $j \in [n]$, $\beta\in \RR$;
\item $z_1^{d_1}\cdots z_n^{d_n}
f(\lambda_1 z_1^{-1},\ldots,\lambda_n z_n^{-1})\in\HH_n(\RR)$ if 
$\pm (\lambda_1,\ldots,\lambda_n)\in\RR_+^n$;
\item $f(z_1,\ldots, z_{i-1}, z_j, z_{i+1}, \ldots, z_n) \in \HH_{n-1}(\RR)$ 
for 
$i\neq j \in [n]$.
\end{itemize}
\item If $\{f_j\}_{j=1}^\infty \subset \HH_n(\RR)$ and 
$f \in \RR[z_1,\ldots, z_n]\setminus \{0\}$ is the limit, uniformly on 
compact subsets of $\CC^n$, of the sequence $\{f_j\}_{j=1}^\infty$ then 
$f \in \HH_n(\RR)$. 
\end{enumerate}
\end{proposition}

Important examples of real stable 
polynomials are given by the following construction 
(see, e.g., \cite{BBS1,BBS2}).

\begin{proposition}\label{det-pen}
Let $A_1,\ldots, A_m$ be (complex) positive semi-definite matrices and let $B$
be a (complex) Hermitian matrix, all matrices being of the same size.
\begin{enumerate} 
\item The polynomial 
\begin{equation}\label{pencil}
z=(z_1,\ldots,z_n)\mapsto f(z)= \det(z_1A_1+\cdots +z_mA_m+B) 
\end{equation}
is either identically zero or real stable; 
\item If $B$ is also 
positive semi-definite then $f$ has all non-negative coefficients.
\end{enumerate}
In particular, if $Z=\text{{\em diag}}(z_1,\ldots,z_m)$ and $A$ is a positive 
semi-definite $m\times m$ matrix 
then $\det(A+Z)$ is a multi-affine real 
stable polynomial with all non-negative coefficients, hence a (positive)
constant multiple of a strongly 
Rayleigh polynomial (cf.~Definition~\ref{def-rsm}). 
\end{proposition}

\begin{proof}
By the second part of Proposition~\ref{pro-rs} and a density argument,  
it suffices to prove (1) 
in the case when  all matrices $A_j$, $1\le j\le m$, are 
positive definite. Set $z(t)=\lambda t+ \mu$ with
$\lambda=(\lambda_1,\ldots,\lambda_n)\in \RR_+^n$,
$\mu=(\mu_1,\ldots,\mu_n)\in \RR^n$, and $t\in \CC$. Note that 
$P:=\sum_{j=1}^m\lambda_jA_j$
is positive definite and thus it is invertible and has a square root $Q$ (recall
that $\RR_+=(0,\infty)$). 
Then 
\begin{equation*}
f(z(t))=\det(P)\det(tI+QHQ^*),
\end{equation*}
where $H:=B+\sum_{j=1}^m\mu_jA_j$ is a Hermitian matrix. 
Therefore, $f(z(t))$ is a polynomial in $t$ that is a 
constant multiple of the characteristic polynomial of a
Hermitian matrix and so it must have all real zeros. By 
the first part of Proposition~\ref{pro-rs}, this proves (1), and then (2)
follows by noticing that the coefficients of the polynomial in \eqref{pencil}
are products of non-negative numbers, namely principal minors of the 
positive semi-definite matrices $A_1,\ldots,A_m,B$.
\end{proof}

Further properties of (not necessarily multi-affine) real stable 
polynomials, their linear
preservers and applications to entire function theory in one or several 
variables, matrix theory, combinatorics, convex optimization and statistical 
mechanics may be found in \cite{BBS,BBS1,BBS2,PB,COSW,guler,LS}.

\subsection{Hadamard-Fischer Inequalities and $\GKK$-Matrices}\label{ss-gkk}

Given an $n\times n$ matrix $A$ over $\CC$ and $S\subseteq [n]$ denote
by $A[S]$ the principal submatrix of $A$ with rows
and columns indexed by $S$ and let $A\langle S\rangle=\det(A[S'])$ be the 
principal minor of $A$ with rows
and columns indexed by $S'=[n]\setminus S$, where 
$A\langle [n]\rangle=\det(A[\emptyset]):=1$. The matrix $A$ is called  a
$P$-{\em matrix} if all its principal minors are positive. A much studied 
class of $P$-matrices that unifies several important types of matrices is 
the class of $\GKK$-matrices, see \cite{FJ,HS} and references 
therein. 
A $P$-matrix $A$ is a $\GKK$-{\em matrix} (after Gantmacher-Krein and
Kotelyansky) if it satisfies the {\em
Hadamard-Fischer-Kotelyansky inequalities}, that is,
\begin{equation}\label{eq-HFK}
A\langle S\rangle A\langle T\rangle \ge A\langle S\cup T\rangle A\langle
S\cap T\rangle,\quad S,T\subseteq [n].
\end{equation}
In other words, an $n \times n$ $P$-matrix $A$  is $\GKK$ if and only if
the probability measure $\mu_A$ on $2^{[n]}$ defined by
\begin{equation}\label{meas-matrix} 
\mu_A(S)= A\langle S \rangle\det(A+I)^{-1}
\end{equation}
satisfies the negative lattice condition (NLC). It is not hard to see that the
generating polynomial of such a measure is
$$
g_{\mu_A}(z) = \det(A+I)^{-1}\sum_{S \subseteq [n]} A\langle S \rangle
z^S=\det(A+I)^{-1}\det(A+Z),\quad z^S:= \prod_{i \in S}z_i,
$$
where  $Z=\text{diag}(z_1,\ldots,z_n)$. Examples of $\GKK$-matrices are:
\begin{itemize}
\item[(I)] Positive definite matrices;
\item[(II)] Totally positive matrices, i.e., matrices for which {\em all}
minors are positive;
\item[(III)] Non-singular $M$-matrices; recall that an $M$-{\em matrix}
(named after Minkowski) is a matrix with all non-negative principal minors 
and all non-positive off-diagonal entries.
\end{itemize}

An ambitious project would be to classify $\GKK$-matrices according to
the various negative dependence properties (introduced in \S \ref{s2-plet}) 
that these matrices satisfy. This would take us beyond the scope of this 
paper. In what follows we will only briefly discuss some of these aspects.

\begin{definition}\label{d-mat-ral}
We say that a $P$-matrix $A$ is a {\em Rayleigh matrix} if the associated 
probability measure $\mu_A$ defined in \eqref{meas-matrix} is Rayleigh 
(in the sense of Definition~\ref{d1}).
\end{definition} 

By \cite[Theorem 4.4]{W}, 
Rayleigh matrices are $\GKK$. In fact, the following holds.

\begin{theorem}\label{R-matrix}
Let $A$ be an $n \times n$ matrix over $\CC$. The following are equivalent:
\begin{itemize}
\item[(1)] $A$ is a Rayleigh matrix;
\item[(2)] $A+X$ is a $\GKK$-matrix for all $X=\diag(x_1,\ldots, x_n)$,
where $x_i \geq 0$ for all $1\leq i \leq n$.
\end{itemize}
\end{theorem}

\begin{proof}
(1) $\Rightarrow$ (2): As we already noted, by \cite[Theorem~4.4]{W} 
Rayleigh matrices are $\GKK$ and by definition the set of all Rayleigh 
matrices is obviously closed under adding positive diagonal matrices.

(2) $\Rightarrow$ (1): Let $f=\det(A+Z)$. Then
\begin{eqnarray*}
&& \frac {\partial f} {\partial z_i}(x) \frac {\partial f} {\partial z_j}(x)
- \frac {\partial^2 f} {\partial z_i \partial z_j}(x) f(x) =\\
&& (A+X)\langle i \rangle \cdot (A+X)\langle j \rangle - (A+X)\langle i,j
\rangle \cdot (A+X)\langle\emptyset \rangle.
\end{eqnarray*}
Using the Hadamard-Fischer-Kotelyansky inequalities \eqref{eq-HFK} for $A+X$ 
with $S=\{i\}$ and $T=\{j\}$ one gets the desired conclusion.
\end{proof}

It follows immediately from  Theorem~\ref{R-matrix} that non-singular
$M$-matrices and positive definite matrices are Rayleigh. Now 
for any $s,t\in(0,1)$ one can easily check that the matrix 
\begin{equation}\label{b-mat}
B=
 \left[\begin{array}{ccc}
 1 & 1/2 & s/4 \\
 1/2 & 1 & 1/2 \\
 t/4 & 1/2 & 1 \\ \end{array} \right]
\end{equation}
 is totally positive.
Given an $n\times n$ matrix $A$ and two subsets $S, T \subseteq [n]$ of the
same size we let $A(S,T)$ denote the minor of $A$ lying in rows indexed by 
$S$ and columns indexed by $T$. It was proved by Gantmacher-Krein \cite{GK} 
and Carlson \cite{Ca}
that a  necessary and sufficient condition for a $P$-matrix to be $\GKK$
is that
$$
A(S,T)A(T,S) \geq 0\text{ for any }S,T \subseteq [n]\text{ with } 
|S|=|T|=|S\cup T|-1.
$$
Since
$$
(B+I)(12,23)= ( 1-2s)/4\text{ and } (B+I)(23,12)=(1-2t)/4
$$
we may choose $s,t\in (0,1)$ so that 
$(B+I)(12,23)(B+I)(23,12)<0$. We conclude that
for such values of $s,t$ the matrix  $B+I$ fails to be $\GKK$.

It is easy to see that the measures associated with positive definite and
totally positive matrices have $\ULC$ rank sequences. This follows from
Newton's inequalities and the fact that 
Hermitian matrices and totally positive matrices have real-rooted
characteristic polynomials (in the case of totally positive matrices this
goes back to Gantmacher and Krein \cite{GK}). For $M$-matrices it is not
as obvious but nonetheless true.
In \cite{HS} Holtz and Schneider conjectured that $M$-matrices satisfy the
conclusion of Conjecture~\ref{big}. This was subsequently confirmed by Holtz
in \cite{HM} using a certain monotonicity property for so-called
symmetrized Fischer 
products of $M$-matrices due to James-Johnson-Pierce \cite{JJP} (see 
\cite{BBS1} for further results on Fischer products).

We end this section with some open questions.

\begin{question}
Is $\mu_A$ an $\ULC$ measure whenever $A$ is a Rayleigh matrix?
\end{question}

\begin{question}
For which $\GKK$-matrices $A$ is $\mu_A$ negatively associated?
\end{question}

Again, negative association fails for measures associated with 
totally positive matrices. Indeed,
let $B$ be the totally positive $3\times 3$ matrix defined in \eqref{b-mat} 
and fix $s,t\in (0,1)$ such that $A:=B+I$ is not $\GKK$ (see the above 
discussion). 
It is known \cite{Ca,GK} that an $n \times n$
$P$-matrix $C$ is $\GKK$ if and only if
\begin{equation}\label{w-kot}
C\langle S \cup\{i\} \rangle \cdot C\langle S \cup\{j\} \rangle \geq
C\langle S \cup\{i,j\} \rangle \cdot C\langle S  \rangle
\end{equation}
for all $S \subseteq [n]$ and $i,j\in [n]\setminus S$ with $i\neq j$.  
Hence, since
$A$ fails to be $\GKK$ it follows that
\eqref{w-kot} fails for $A$ and some $S, i,j$. Unraveling the definitions 
this means that for some principal submatrix $B'$ of
$B$ the corresponding measure $\mu_{B'}$ will fail to have pairwise negatively
correlated variables (cf.~Definition~\ref{d-pnc}). In particular, since
NA $\Rightarrow$ p-NC, the measure $\mu_{B'}$ cannot be negatively 
associated. 

\begin{question}
Describe all $n\times n$ Rayleigh  matrices.
\end{question}

\begin{question}
Characterize the class of all real stable $n\times n$
matrices, that is,
$n\times n$ matrices $A$ such that $\det(A+Z)$ is real stable.
\end{question}

\subsection{Determinantal Probability Measures}\label{ss-det-m}

In \cite{Lyons} Lyons defined a {\em determinantal probability measure} 
on $2^{[n]}$ to 
be a measure $\mu\in\fP_n$ such that
there is an $n \times n$ matrix $A$ so that for any subset $S$ of 
$[n]$ one has 
\begin{equation}\label{det-me}
\mu(\{T: S \subseteq T\}) = \det(A[S]),
\end{equation}
where as before $A[S]$ denotes the principal submatrix of $A$ whose rows and 
columns are indexed by $S$. As discussed in e.g.~\cite{BKPV,Lyons}, such 
measures have important applications ranging from 
fermionic processes/continuous scaling limits of discrete point 
processes \cite{DVJ,Lyons-Steif} to the distribution of eigenvalues of 
random matrices and zeros of the Riemann zeta function \cite{Con}, 
transfer current matrices \cite{BP}, 
non-intersecting random walks \cite{Joh2}, 
Poissonized versions of the Plancherel
measure on partitions \cite{BOO}
and Young diagrams \cite{OR}, discrete orthogonal polynomial 
ensembles \cite{Joh}, etc.

Recall that a square matrix $A$ is said to be a {\em contraction} 
if $||A||\le 1$, where $||\cdot||$ denotes the supremum (operator) norm. 
Clearly, a positive semi-definite matrix is a contraction if and only
if all its eigenvalues are in the interval $[0,1]$. Such matrices were
called {\em positive contractions} in \cite{Lyons}. 

\begin{theorem}\label{contract}
Suppose that  $\mu$ is a determinantal measure on $2^{[n]}$ whose
corresponding matrix is a positive contraction. Then $\mu$ is $\CNA+$.
\end{theorem}

This is one of the main results of \cite{Lyons} 
(cf.~\cite[Theorem 8.1]{Lyons}), where 
exterior (Grassmann) algebra methods were used to investigate determinantal 
probability measures and, in particular, to prove the above theorem.
As noted in {\em loc.~cit.}, Theorem~\ref{contract} provides a powerful tool 
for studying determinantal probability measures, comparable to the FKG theorem
in the theory of positive association, and is crucial to most of the results
in e.g.~\cite{Lyons-Steif}. 

Let us show that Theorem \ref{contract} actually follows from 
Theorem \ref{stab-NA} of \S \ref{s-sna} and the fact that determinantal 
measures induced by positive contractions are strongly Rayleigh: 

\begin{proposition}\label{det-stab}
Suppose that  $\mu$ is a determinantal measure on $2^{[n]}$ whose
corresponding matrix is a positive contraction. Then $\mu$ is strongly 
Rayleigh.
\end{proposition}

\begin{proof}
Since positive definite matrices are dense in the set of all positive 
semi-definite matrices, by Proposition \ref{pro-rs} (3) it is enough to prove 
the statement for invertible
positive contractions. Elementary computations show that 
$$g_\mu(z_1,\ldots,z_n)=\det(I-A+AZ)=\det(A)\cdot \det(A^{-1}-I+Z),$$ 
where $Z=\text{diag}(z_1,\ldots,z_n)$. 
Since $A$ is an invertible positive contraction, $A^{-1}-I$ is positive 
semi-definite and
$\det(A)>0$, so by \eqref{pencil} the polynomial $g_\mu$ is 
real stable. 
\end{proof}

Another main result of \cite{Lyons} is that if $A$ and $B$ are commuting 
positive contractions and $A \leq B$, 
then the measure corresponding to $B$ stochastically dominates
the one corresponding to $A$ (Definition \ref{d-stoc}). In 
\S \ref{ss-trunc} we extend this result using the theory of stability and 
strongly Rayleigh measures and we prove that the condition 
$[A,B]=0$ (where $[\cdot,\cdot]$ denotes the Lie bracket)
can be dropped from the hypothesis. 

\subsection{Graphs, Laplacians, Spanning Trees and the Random Cluster 
Model}\label{ss-trees}

Let $G=(V,E)$ be a graph with vertex set $V=[n]$ and edge set $E$. Associate 
to each edge $e \in E$ a variable $w_e$. If $e$ connects $i$ and $j$ let 
$A_e$ be the $n \times n$ positive semi-definite matrix with the $ii$-entry 
and 
the $jj$-entry equal to $1$, the $ij$-entry and $ji$-entry equal to $-1$, 
and all other entries equal to $0$. The {\em Laplacian}, $L(G)$,  of $G$ is 
defined by 
$$
L(G)= \sum_{e \in E}w_eA_e.
$$
Let
$$
f_G(z,w)= \det(L(G)+Z).
$$
Thus, by \eqref{pencil}, $f_G$ is a multi-affine real stable polynomial with 
non-negative coefficients.   
The Principal Minors Matrix-Tree Theorem (see, e.g., \cite{Chaiken}) says that 
$$
 f_G(z,w)= \sum_{F}z^{\roots(F)}w^{\edges(F)}, 
$$ 
where the sum is over all rooted spanning forests $F$ in $G$, $\roots(F) 
\subseteq V$ is the set of roots 
of $F$ and $\edges(F) \subseteq E$ is the set of edges used in $F$. 
Alternatively, we may write $f_G$ as 
$$
f_G(z,w)= \sum_{F} w^F \prod_{C \in \mathcal{C}(F)} (\sum_{j \in C}z_j),
$$
where $\mathcal{C}(F)$ is the partition of $V$ whose parts are the vertices 
in the maximal trees in $F$ and the sum is over all forests in $G$.  Since 
the class of stable polynomials is closed under differentiation and 
specialization of  variables at real values (see, e.g., \cite{BBS2}) we have 
that the {\em spanning tree polynomial} $T_G(w)= \sum_Tw^T$, where the sum is 
over all spanning trees, is real stable (which is widely known).  The 
corresponding probability measure is usually called the uniform random 
spanning tree measure and is consequently strongly Rayleigh. Now if 
$\ben$ denotes the vector of all ones as in \S \ref{s2-plet} then 
$$
f_G(\ben,w)= \sum_{F} \omega_G(F)w^F,
$$
where the sum is over all spanning forests (independent sets) in $G$ and 
 $\omega_G(F)$ is the product of the sizes of the connected components of 
$F$. In \cite[Conjecture~5.10]{W} Wagner conjectured that $f_G(\ben,w)$ is a 
Rayleigh polynomial. By the above discussion we have shown that more is 
actually true, namely:

\begin{theorem}\label{t-fg}
For any graph $G$ the polynomial $f_G(\ben,w)$ is real stable.
\end{theorem}

Let $G=(V,E)$ be a graph and let $q>0$ and $0 \leq p \leq 1$ be parameters. 
The {\em random cluster (RC) model}, see \cite{grimmett},  
is the measure $\mu$ on $2^E$ defined 
by 
$$
\mu(F) = \frac 1 {Z_{\text{RC}}} q^{k(F)}p^{|F|}(1-p)^{|E \setminus F|},  
\quad F \subseteq E, 
$$
where $k(F)$ is the number of connected components in the graph $(V,F)$ and 
$Z_{\text{RC}}$ is the appropriate normalizing factor.  
The uniform spanning tree measure can be obtained
from this model by letting $p\to 0$ and $q/p \to 0$, see 
\cite[Section 1.5]{grimmett}. If
$q\geq 1$, the
RC model satisfies PLC and is therefore positively associated. 
If $q\leq1$, NLC is
satisfied, but other negative correlation properties are largely a matter of
conjecture. 
We are interested in properties closed under external fields so we assume 
that $p=1/2$. The generating polynomial of $\mu$ is then a constant multiple 
of the multivariate Tutte polynomial (see \cite{sok}):
$$
Z_G(z,q)= \sum_{F \subseteq E}q^{k(F)}z^F, \quad z=(z_e)_{e\in E}.
$$

If $G=C_n, n \geq 3$, is a cycle and $E=[n]$, then  
$$
Z_G(z,q)= \prod_{j=1}^n(q+z_j) +(q-1)z_1\cdots z_n 
$$
and 
$$
\frac{\partial Z_G}{\partial z_1}\frac{\partial Z_G}{\partial
z_2}-Z_G\frac{\partial^2 Z_G}
{\partial z_1\partial z_2}=  q^2(1-q)\prod_{j=3}^n z_j(q+z_j),
$$
so in this case the RC measure is Rayleigh, but not strongly Rayleigh. 

If $G$ is a tree, then 
$$
Z_G(z,q)= q\prod_{e\in E}(q+z_e),
$$
so the RC model is strongly Rayleigh. 
We conclude that on a
general graph $G$, the RC model is strongly Rayleigh if and only if $G$ is 
acyclic.  

With the help of {\em Mathematica}$^\copyright$, the
third author has verified that the RC model is Rayleigh for all simple 
graphs with at most five vertices. 

\subsection{Interacting Particle Systems and Exclusion 
Evolutions}\label{ss-ex-case}

The exclusion process is one of the main models considered in the area of
Probability Theory known as Interacting Particle Systems. The idea is 
that particles move in continuous time on a countable set $S$, in such a way 
that there is always at most one particle per site. The motion of
each particle would be a continuous time Markov chain on $S$, except that
transitions to occupied sites are forbidden. This process has been used to 
model many situations, including biopolymers and traffic flow. We refer to
\cite[Chap.~VIII]{L0} and \cite[Part 3]{L00} for detailed
treatments of this process. 

The (finite) {\em symmetric} exclusion process
generates a continuous time evolution on the set $\fP_n$. 
The limiting measure under the 
symmetric exclusion evolution with initial measure $\mu$ 
is the symmetrization $\mu_s$ of $\mu$ as defined in \S \ref{s21}. 
A problem that has attracted some attention is the following (see, e.g.,
\cite{L0,L2}).

\begin{problem}\label{p-see}
Find a natural negative dependence property that is preserved by symmetric 
exclusion evolutions.
\end{problem}

Since the limiting distribution of the evolution as $t\to\infty$ is the 
symmetrization
of the initial distribution defined in \eqref{sym}, any such property would  
have to be preserved by
this symmetrization procedure. In order for a solution to Problem~\ref{p-see}
to be useful in settling
Conjecture~\ref{con-lig}, the negative dependence property should be 
satisfied  by non-random/deterministic distributions and should
imply NA. Problem~\ref{p-see} is also
motivated by a theorem of Harris \cite[Theorem 2.14]{L0} that provides a 
general criterion for the preservation of positive dependence for Markov 
processes. Theorems 3.1 and 3.3 of \cite{L2} give negative dependence 
properties that are preserved by the symmetric exclusion evolution. However, 
one can make an argument that neither of these is ideal -- the first is too 
weak while 
the second is too strong. In \cite[Theorem 3.2]{L2} it was proved that NA is
not preserved by the symmetric exclusion 
process. The examples that we construct in \S \ref{s-countex}
further show that in fact none of the properties NLC, h-NLC, Rayleigh/h-NLC+, 
CNA, CNA+ defined in \S \ref{s21} is preserved by such evolutions. 
This is a consequence of the 
following criterion established by
Pemantle in \cite[Theorem 2.7]{P} (see also \cite[Proposition 3.6]{W}) for the 
ultra log-concavity of rank sequences of 
symmetric polynomials in $\PP_n$:

\begin{theorem}[\cite{P}]\label{pro-w}
For $0\le k\le n$ let $e_k(z_1,\ldots,z_n)$ be the $k$-th elementary symmetric 
function on $\{z_1,\ldots,z_n\}$ and consider the polynomial
$$f(z_1,\ldots,z_n)=\sum_{k=0}^{n}a_ke_k(z_1,\ldots,z_n),$$
where $a_k\ge 0$, $0\le k\le n$. The following are equivalent:
\begin{enumerate}
\item $f$ satisfies the conclusion of 
Conjecture~\ref{big}, that is, it has an ULC rank sequence 
(cf.~Definition~\ref{d-ulc});
\item $f$ has either (and then all) of the following five properties: NLC, 
h-NLC, Rayleigh/h-NLC+, CNA, CNA+.
\end{enumerate}
\end{theorem}

In \S \ref{n-s6} we show that the strongly Rayleigh 
property is preserved by (finite) symmetric exclusion evolutions  
and thus provide an answer to Problem~\ref{p-see}. 
We then prove that the distribution at time 
$t>0$ for the (not necessarily finite) symmetric exclusion process is 
CNA+ whenever the initial distribution is a product measure, which confirms 
Conjecture \ref{con-lig}.

Symmetric strongly Rayleigh measures are described in the following theorem. 
Recall that a (possibly infinite) matrix 
$A$ is said to be totally 
non-negative of order $p$ or $\text{TP}_p$ for short 
if for any $k\le p$ all its $k \times k$ minors are non-negative \cite{Kar}.
If $A$ has all non-negative minors then $A$ is called 
totally non-negative or $\TP$ for short. 

\begin{theorem}[Exchangeable Strongly Rayleigh Case]\label{ESC}
Suppose $\mu\in\PP_n$ has a symmetric generating polynomial and 
rank sequence $\{r_k\}_{k=0}^n$. The following are equivalent:
\begin{itemize}
\item[(a)] The generating polynomial of $\mu$ is stable, i.e., $\mu$ is
strongly Rayleigh; 
\item[(b)] All zeros of the univariate polynomial $\sum_{k=0}^{n}r_kz^k$
are real; 
\item[(c)] The infinite Toeplitz matrix $(r_{i-j})_{i,j=0}^{\infty}$ is $\TP$,
where $r_k=0$ if $k\notin \{0\}\cup [n]$. 
\end{itemize}
\end{theorem}

\begin{proof}
The equivalence (b) $\Leftrightarrow$ (c) is a classical result due to 
Aissen, Schoenberg and 
Whitney \cite{ASW}, while the equivalence (a) $\Leftrightarrow$ (b) is a 
consequence of the Grace-Walsh-Szeg\"o coincidence theorem 
(see Theorem~\ref{GWS} below). 
\end{proof}

\section{Negative Dependence Theory for Strongly Rayleigh 
Measures}\label{s-stable}

The importance of strongly Rayleigh measures -- that is, measures with stable 
generating polynomials (cf.~Definition~\ref{def-rsm}) -- in negative 
dependence stems from the theorem below proved by one of the authors  
in \cite{PB} (cf.~\S \ref{ss-real-s}).  

\begin{theorem}\label{stable-Rayleigh}
Let $g \in \RR[z_1,\ldots, z_n]$ be a multi-affine polynomial. Then $g$ is 
stable if and only if 
$$
\frac {\partial g} {\partial z_i}(x) \frac {\partial g} {\partial z_j}(x) 
\geq \frac {\partial^2 g} {\partial z_i \partial z_j}(x) g(x)$$
for any $x=(x_1,\ldots,x_n)\in \RR^n$ and $i,j\in [n]$.
\end{theorem}

Hence, perhaps somewhat surprisingly, if the Rayleigh condition 
in Definition~\ref{d1} is enforced to hold for all real vectors then this 
stronger correlation condition is equivalent to real stability. 
Clearly, a strongly Rayleigh measure is automatically Rayleigh 
(cf.~Remark~\ref{r-term}). As we will see in this section, the 
geometric properties encoded in the concept of real stability allow us 
to show that strongly Rayleigh measures 
enjoy all virtues of negative dependence, including the strongest
form of negative association $\CNA+$. 

\subsection{Symmetric Homogenization via G\aa rding Hyperbolic 
Polynomials}\label{s-gar}

The following theorem shows that the symmetric 
homogenization (cf.~\S \ref{s22}) of any strongly Rayleigh measure is again
strongly Rayleigh, so that 
strongly Rayleigh measures belong to the class PHR (Definition~\ref{phr}).
 
\begin{theorem}\label{SH}
If $\mu\in \fP_{n}$ is strongly Rayleigh then so is its symmetric 
homogenization $\mu_{sh}\in\fP_{2n}$.
\end{theorem}

The proof of Theorem~\ref{SH} makes use of the theory of hyperbolic 
polynomials that has its origins in partial differential equations 
and was developed by 
Petrovsky, G\aa rding, H\"ormander, Atiyah and Bott 
\cite{ABG,garding,hormander}. Recently, hyperbolic polynomials have proved
to be quite useful in other areas of mathematics such as convex optimization 
\cite{BGLS,guler}, complex analysis \cite{BBS}, matrix theory and 
combinatorics \cite{Gu1,Gu}. 

Let $e\in \RR^n$. A homogeneous polynomial $p \in \RR[z_1,\ldots, z_n]$ is 
said to be {\em (G\aa rding) hyperbolic with respect to} $e$ if $p(e)\neq 0$ 
and for all $x \in \RR^n$ the univariate polynomial 
$t \mapsto p(x+te)$ has only real zeros. Recall that the {\em homogenization} 
of a polynomial $f(z_1,\ldots,z_n)=\sum_{\alpha \in \NN^n}a(\alpha)z^\alpha  
\in \CC[z_1, \ldots, z_n]$ of degree $d$ is given by
 $$f_H(z_1,\ldots,z_{n+1}) = z_{n+1}^df(z_1/z_{n+1}, \ldots, z_n/z_{n+1}).$$  

Recall the characterization of real stable polynomials given in 
Proposition~\ref{pro-rs} (1). The relationship between 
real stable polynomials and hyperbolic polynomials is made explicit by the 
following result, see \cite[Proposition 1]{BBS2}.

\begin{proposition}\label{hyp-stab}
Let $f \in \RR[z_1,\ldots, z_n]$. Then $f$ is real stable if and only if 
$f_H$ is hyperbolic with respect to all vectors $e \in \RR^{n+1}$ of the form 
$e=(e_1, \ldots, e_n, 0)$, where $e_i>0$ for $1 \leq i \leq n$. 
\end{proposition} 

Let $p$ be hyperbolic with respect to $e\in\RR^n$. The cone 
$$
C_e(p)=\{x \in \RR^n : p(x+te) \neq 0, t \geq 0\}
$$
is called the {\em hyperbolicity cone} of $p$. The following fundamental 
properties of the hyperbolicity cone are due to G\aa rding \cite{garding}, 
see also \cite{ABG,guler,hormander}. 

\begin{proposition}\label{fund-hyp}
Let $p \in \RR[z_1,\ldots,z_n]$ be hyperbolic with respect to $e\in\RR^n$. 
Then 
\begin{itemize}
\item[(a)] $C_e(p)$ is convex; 
\item[(b)] $C_e(p)$ is equal to the connected component of the set 
$\{ x \in \RR^n : p(x) \neq 0\}$ that contains $e$;
\item[(c)] $(s,t)\mapsto p(x+su+te)$ is real stable for any $x\in \RR^n$ 
and $u\in C_e(p)$; 
\item[(d)] $p$ is hyperbolic with respect to any $u\in C_e(p)$ and 
$C_u(p)=C_e(p)$.
\end{itemize}
\end{proposition}

Below are examples of hyperbolic polynomials and their hyperbolicity 
cones: 

\begin{itemize}
\item[(I)] The polynomial $p(z_1,\ldots,z_n)= z_1z_2\cdots z_n$ is hyperbolic 
with respect to the direction $e=(1,\ldots,1)\in\RR^n$ and the hyperbolicity 
cone is the positive orthant $C_e(p)=\RR_+^n$; 
\item[(II)] The (principal) symbol of the hyperbolic wave equation, i.e., 
the polynomial 
$$
p(z_1,\ldots, z_n) = z_1^2-\sum_{j=2}^n z_j^2
$$
is hyperbolic with respect to the direction $e=(1,0,\ldots,0)\in\RR^n$ and 
the hyperbolicity cone is the Lorentz cone 
$$
C_e(p)=\left\{z \in \RR^n : \sqrt{z_2^2+\cdots+z_n^2} \leq z_1\right\};
$$
\item[(III)] Let $z_{ij}$, $1\le i\le j\le n$, be $n(n+1)/2$ different 
variables and let $X=(x_{ij})_{i,j=1}^n$  be the matrix with entries 
$$
x_{ij} = \begin{cases} 
z_{ij} \mbox{ if } i \leq j, \\
z_{ji} \mbox{ if } i > j.
\end{cases} 
$$
The polynomial (in $n(n+1)/2$ variables) 
$p(z_{11}, \ldots, z_{nn}) = \det(X)$ is  hyperbolic with respect to the
$n\times n$ identity matrix $I$ and the hyperbolicity cone is the cone of 
positive definite $n\times n$ matrices; 
\item[(IV)] Let $A_1, \ldots, A_m$ be symmetric $n \times n$ matrices and let 
$e=(e_1,\ldots,e_m) \in \RR^m$ be such that  $e_1A_1 + \cdots+e_mA_m$ is 
positive definite. Then the polynomial 
$$
p(z_1,\ldots, z_m) = \det\left(\sum_{j=1}^m z_jA_j\right)
$$
is hyperbolic with respect to $e$ and the hyperbolicity cone is given by 
$$
C_e(p)= \left\{ x \in \RR^m : \sum_{j=1}^n x_jA_j \mbox{ is positive 
definite}\right\}.
$$   
\end{itemize}

\begin{theorem}\label{homhom}
Suppose that all the coefficients of $f \in \RR[z_1,\ldots, z_n]$ are 
non-negative. The following are equivalent:
\begin{enumerate} 
\item $f_H$ is stable; 
\item $f$ is stable; 
\item $f_H$ is hyperbolic with respect to some vector $e$ with $e_i\geq 0$; 
$1 \leq i \leq n+1$;
\item $f_H$ is hyperbolic with respect to any vector $e$ with $e_i>0$, 
$1 \leq i \leq n+1$. 
\end{enumerate}
\end{theorem}

\begin{proof}
By the fact that the stability property is 
closed under setting variables equal to real numbers 
(Proposition~\ref{pro-rs}) and Proposition~\ref{hyp-stab} we have that 
(1) $\Rightarrow$ (2) $\Rightarrow$ (3). 
Now if (3) holds then since $f$ has all non-negative coefficients it follows
from Proposition~\ref{fund-hyp} (b) that 
$C_e(f_H)$ contains the cone $\RR_+^{n+1}$, which by 
Proposition~\ref{fund-hyp} (d) proves (4). Finally, 
(4) $\Rightarrow$ (1) by Proposition~\ref{pro-rs}.
\end{proof}

Let $f \in \CC[z_1,\ldots, z_n]$ be a polynomial of degree $d_i$ in $z_i$, 
$1\le i\le n$. The 
{\em polarization} $\pi(f)$ of $f$ is the unique polynomial in 
$\sum_{i=1}^{n}d_i$ variables 
$z_{ij}$, $1 \leq i \leq n$, $1 \leq j \leq d_i$, satisfying:
\begin{enumerate}
\item[(i)] $\pi(f)$ is multi-affine; 
\item[(ii)] $\pi(f)$ is symmetric in the variables 
$z_{i1}, \ldots z_{id_i}$ for 
$1 \leq i \leq n$;
\item[(iii)] if we let  $z_{ij}=z_i$ for all $i,j$ in $\pi(f)$ we recover $f$.
\end{enumerate}

Recall that a {\em circular region} in $\CC$ is either an open or closed 
affine half-plane or the open or closed interior or exterior of a circle. 
The famous Grace-Walsh-Szeg\"o coincidence theorem 
\cite{grace,RS,walsh} is stated next. For a proof of the following version we 
refer to \cite[Theorem 2.12]{COSW}. 

\begin{theorem}[Grace-Walsh-Szeg\"o]\label{GWS}
Let $f \in \CC[z_1,\ldots, z_n]$ be symmetric and multi-affine and let $C$ be 
a circular region 
containing the points $\zeta_1, \ldots, \zeta_n$. Suppose that either $f$
has total degree $n$ or $C$ is convex (or both). Then there exists at least
one point 
$\zeta \in C$ such that 
$f(\zeta_1, \ldots, \zeta_n)=f(\zeta, \ldots, \zeta)$.

\end{theorem}
From the Grace-Walsh-Szeg\"o coincidence theorem we deduce:

\begin{corollary}\label{polarizeit}
Let $f \in \CC[z_1,\ldots, z_n]$. Then 
$f$ is stable if and only if $\pi(f)$ is stable. 
\end{corollary}

\begin{proof}
If $\pi(f)$ is stable then so is $f$ since the latter is recovered from 
$\pi(f)$ by setting variables equal as in (iii) above. Now suppose that 
$\pi(f)$ is not stable. Then there are numbers 
$\zeta_{ij} \in \{z \in \CC : \Im(z)>0\}$, $1 \leq i \leq n$, 
$1 \leq j \leq d_i$, such 
that $\pi(f)(\zeta)=0$, where
$\zeta
=(\zeta_{11},\ldots,\zeta_{1d_1},\ldots,\zeta_{n1},\ldots,\zeta_{nd_n})$.
By successively applying the Grace-Walsh-Szeg\"o coincidence theorem we see 
that there is a vector $\zeta' = 
(\zeta_{11}',\ldots,\zeta_{1d_1}',\ldots,\zeta_{n1}',\ldots,\zeta_{nd_n}')$ 
with coordinates in $\{z \in \CC : \Im(z)>0\}$ for which $\zeta_{ij}'
= \zeta_{k\ell}'$ whenever $i=k$ and $\pi(f)(\zeta')=\pi(f)(\zeta)=0$. Hence 
$f(\zeta_{11}', \zeta_{21}' \ldots, \zeta_{n1}')=\pi(f)(\zeta')=0$, so $f$ is 
not stable either. 
\end{proof}

We now have all the tools to prove Theorem~\ref{SH}. 

\begin{proof}[Proof of Theorem~\ref{SH}]
Let $g$ be the generating polynomial of $\mu\in\fP_n$ and 
suppose that $g$ is stable of degree $d$. 
Then so is the homogenization $g_H$ of $g$ by Proposition~\ref{hyp-stab}. 
Now straightforward computations show that the generating polynomial of 
the symmetric homogenization $\mu_{sh}\in\fP_{2n}$ is given by
 $\pi(z_{n+1}^{n-d}g_H)$, so the theorem follows from 
Corollary~\ref{polarizeit}. 
\end{proof} 

\subsection{Stability Implies Negative Association}\label{s-sna}

The main argument in the proof of the following theorem goes back to 
Feder-Mihail \cite{feder} although the proof of the full theorem is scattered 
over the literature \cite{DJR,Lyons,Lyons-Peres,W}. 

\begin{theorem}\label{FM}
Let $\mathcal{S}$ be a class of probability measures satisfying: 
\begin{enumerate}
\item Each $\mu \in \mathcal{S}$ is a measure on $2^{E}$, where $E$ is a 
finite subset of $\{1,2,\ldots\}$ depending on $\mu$; 
\item $\mathcal{S}$ is closed under conditioning; 
\item For each $\mu \in \mathcal{S}$ the variables $X_e$, $e \in E$, 
are pairwise negatively correlated; 
\item Each $\mu \in \mathcal{S}$ has a homogeneous generating polynomial.  
\end{enumerate} 
Then all measures in $\mathcal{S}$ are $\CNA$ (conditionally negatively 
associated).
\end{theorem}

The homogeneity condition in (4) above is essential and ensures 
the existence of
a variable of so-called ``positive influence'' \cite{DJR} (various 
notions of ``influence'' have been studied by Bourgain {\em et al} in 
\cite{BKKKL}). 

A sketch of the proof of Theorem~\ref{FM} is as follows.  
Recall from \S \ref{s21} the definition of an increasing event. Since any 
increasing function depending on a set $E_0$ can be written as a positive 
linear 
combination of characteristic functions, $\chi_{\mathcal{A}}$, of increasing 
events $\mathcal{A}$ depending on $E_0$, it is enough to prove the theorem for 
such characteristic functions. The proof of negative association in the case
when 
$F =\chi_{\mathcal{A}}$ and $G = X_e$ in Definition~\ref{d-na} are such that
$\cA$ does not depend on $e$ is recovered from Feder-Mihail's proof 
\cite{feder,Lyons-Peres}. The general case -- when
$F =\chi_{\mathcal{A}_1}$ and $G = \chi_{\mathcal{A}_2}$ in 
Definition~\ref{d-na} are such that $\cA_1$ and $\cA_2$ depend on disjoint
sets of variables -- then follows from the proof of \cite[Theorem 6.5]{Lyons}. 

Using Theorem~\ref{FM} we can establish the following result.

\begin{theorem}\label{stab-NA}
If $\mu\in \fP_n$ is strongly Rayleigh then it is $\CNA+$. 
\end{theorem}

\begin{proof}
Let $\mathcal{S}$ be the class of all probability measures $\mu$ such that 
\begin{itemize} 
\item each $\mu \in \mathcal{S}$ is a measure on $2^{E}$, where $E$ is a 
finite subset of $\{1,2,\ldots\}$ depending on $\mu$, and 
\item $\mu$ has a stable homogeneous generating polynomial.  
\end{itemize}
Since stable homogeneous polynomials are pairwise negatively correlated and
are also closed under conditioning (cf.~\S \ref{s21}), the class 
$\mathcal{S}$ satisfies all the hypotheses required in Theorem~\ref{FM}.
Therefore, all the measures in $\mathcal{S}$ are negatively associated. Let 
further $\widehat{\mathcal{S}}$ 
be the class of all probability measures $\mu$ such that 
\begin{itemize} 
\item each $\mu \in \widehat{\mathcal{S}}$ is a measure on $2^{E}$, where 
$E$ is a 
finite subset of $\{1,2,\ldots\}$ depending on $\mu$, and 
\item $\mu$ has a stable generating polynomial.  
\end{itemize}
Now by Theorem~\ref{SH} every measure in $\widehat{\mathcal{S}}$ is the 
projection of a measure in 
$\mathcal{S}$. Since negative association is closed under projections and 
the strongly Rayleigh property is closed under conditioning and external 
fields, the theorem follows.
\end{proof}

Note that since the Rayleigh/h-NLC+ property is also closed under conditioning,
projections, and external fields, Theorem~\ref{FM} and slight modifications of
the above arguments actually yield the following stronger version of 
Theorem~\ref{stab-NA}.

\begin{theorem}\label{PHR-NA}
If $\mu\in \fP_n$ is $\PHR$ then it is $\CNA+$. 
\end{theorem}

Since homogeneous Rayleigh polynomials need not be stable, 
strongly Rayleigh $\Rightarrow$ PHR strictly. As we will see in 
\S \ref{s-countex}, one also has PHR $\Rightarrow$ CNA+ strictly.

\subsection{Stochastic Domination and Truncations}\label{ss-trunc}

Recall the notion of stochastic domination from Definition~\ref{d-stoc} in 
\S \ref{s24}. 
It turns out that for strongly Rayleigh measures, stochastic domination is 
intimately connected with central notions in the theory of stable and 
real-rooted polynomials, namely the notions of interlacing zeros and 
proper position \cite{BBS,BBS1,BBS2,PB,Le,RS}. In particular, this connection
allows us to extend one of Lyons' main results \cite{Lyons} and to answer
Pemantle's questions on stochastic domination for truncations of ``negatively
dependent'' measures \cite{P} (cf.~\S \ref{s24}).

\subsubsection{Partial Orders for Strongly Rayleigh Measures and 
Polynomials}\label{sss-stoc}

Let $\alpha_1 \leq \alpha_2 \leq \cdots \leq \alpha_n$ and 
$\beta_1 \leq \beta_2 \leq \cdots \leq \beta_m$ be the zeros of two 
real-rooted polynomials $p,q \in \RR[z]$. These zeros are  
{\em interlacing} if they can 
be ordered so that either   $\alpha_1 \leq \beta_1 \leq \alpha_2 \leq 
\beta_2 \leq 
\cdots$ or $\beta_1 \leq \alpha_1 \leq \beta_2 \leq \alpha_2 \leq 
\cdots$. If 
the zeros of $p,q$ interlace then 
the {\em Wronskian} 
$W[p,q]:=p'q-pq'$ is either non-negative or non-positive on the whole
 real axis 
$\RR$, see, e.g., \cite{Le,RS}. 
We say that $p$ and $q$ 
are  in  
{\em proper position}, denoted $p \ll q$, if the zeros of $p$ and $q$ 
interlace and 
 $W[p,q] \leq 0$. For technical reasons, we also say that the zeros of the 
polynomial $0$ 
interlace the zeros of any (non-zero) real-rooted polynomial and write 
$0 \ll p$ and $p \ll 0$. 

The following theorem proved in \cite{BBS2} 
provides a notion of proper position for  multivariate polynomials. 

\begin{theorem}\label{addz}
Let $f, g \in \RR[z_1,\ldots,z_n]$. The following are equivalent: 
\begin{enumerate}
\item The polynomial $g+if \in \CC[z_1,\ldots,z_n]$ is stable; 
\item The polynomial $g+z_{n+1}f \in \RR[z_1,\ldots,z_n,z_{n+1}]$ is 
real stable; 
\item For all $\lambda \in \RR_+^n$ and $\mu \in \RR^n$ we have 
$$
f(\lambda t + \mu) \ll g(\lambda t + \mu),  
$$
where $f(\lambda t + \mu)
=f(\lambda_1 t+\mu_1, \ldots, \lambda_n t+\mu_n)$, $t\in\CC$.
\end{enumerate}
\end{theorem}

We say that the polynomials $f, g \in \RR[z_1,\ldots,z_n]$ are in 
{\em proper position}, written $f \ll g$, if any of the equivalent 
conditions in Theorem \ref{addz} are satisfied. 

\begin{remark}\label{dada}
The typical proper position is between a real stable polynomial and any of its 
partial derivatives: 
if $f \in \RR[z_1,\ldots,z_n]$ is real stable then 
$\partial_j f \ll f$ for all $j\in [n]$. 
If $f \ll g$ and $\alpha$ is a real number then either 
$f(\alpha,z_2,\ldots,z_n) \ll g(\alpha,z_2,\ldots,z_n)$ or  
$f(\alpha,z_2,\ldots,z_n)= g(\alpha,z_2,\ldots,z_n)\equiv 0$, see \cite{BBS2}.
Note also that by Theorem~\ref{addz} (3) one has $f\ll g\Leftrightarrow 
\alpha f\ll \beta g$ for any $\alpha,\beta\in \RR_+$.
\end{remark} 

We define a partial order $\unlhd$ on the set of strongly Rayleigh 
probability measures on $2^{[n]}$ as follows. 
 
\begin{definition}\label{d-order}
 Let $\SR_n$ be the set of  strongly Rayleigh probability measures on 
$2^{[n]}$. First, define a (preliminary) partial order $\unlhd'$ on $\SR_n$ 
by setting $\mu \unlhd' \nu$ if there exists a sequence of  strongly 
Rayleigh probability measures $\mu=\mu_0, \mu_1, \ldots, 
 \mu_\ell=\nu$ such that 
 $g_{\mu_0} \ll g_{\mu_1} \ll \cdots \ll g_{\mu_\ell}$. Then 
 define the partial order $\unlhd$ on $\SR_n$ as the closure of 
$\unlhd'$, i.e., set $\mu \unlhd \nu$ if there are two sequences 
 $\{\mu_j\}_{j=0}^{\infty}, \{\nu_j\}_{j=0}^{\infty} \subset \SR_n$ such that 
 $$
 \lim_{j \rightarrow \infty} \mu_j =\mu,\, 
 \lim_{j \rightarrow \infty} \nu_j =\nu, \mbox{ and } 
\mu_j \unlhd' \nu_j \text{ for all } j \in \NN 
 $$
(as the measures are defined 
on $2^{[n]}$ it does not matter what notion of limit we use). 
We will also consider $\unlhd$  to be a partial order on 
 the set of generating polynomials of strongly Rayleigh measures and write 
$g_{\mu} \unlhd g_{\nu}$ whenever $\mu \unlhd \nu$.   
\end{definition}

\begin{remark}
It seems likely that the partial orders $\unlhd'$ and $\unlhd$ on $\SR_n$ 
defined above are actually the same. However, since it is 
irrelevant for the present purposes we do not further pursue this 
issue here.
\end{remark}

As we noted in \S \ref{s24}, the stochastic domination relation $\pc$ defines
a partial order on $\fP_n$ (cf.~\cite{MO}). 
The anti-symmetry of $\unlhd$ follows from our next result. 

\begin{proposition}\label{stdom}
Let $\mu$ and $\nu$ be strongly Rayleigh probability measures on $2^{[n]}$. 
Then 
$\mu \unlhd \nu \Longrightarrow \mu \pc \nu$.
\end{proposition}

\begin{proof}
Clearly, we may assume that $g_\mu \ll g_\nu$, 
so that by Theorem~\ref{addz} the
polynomial in $n+1$ variables 
$G=(g_\nu +z_{n+1}g_\mu)/2$ is stable. Let $\gamma\in\fP_{n+1}$ be the 
measure with generating polynomial $G$. By 
Theorem \ref{stab-NA}, $\gamma$ is NA. Let $\AS$ be an 
increasing event not depending on the last coordinate  and 
set $\BS=\{ S  \in 2^{[n+1]} : n+1 \in S \}$. Then 
$$
\frac{1}{2} \mu(\AS) = \gamma(\AS \cap \BS) \leq \gamma(\AS) 
\gamma(\BS)
= \left(\frac{1}{2} \mu(\AS) + \frac{1}{2} \nu(\AS)\right)\frac{1}{2}, 
$$
which gives $\mu(\AS) \leq \nu(\AS)$. 
\end{proof}

Recall that the standard partial order on the set of all Hermitian $n\times n$
matrices -- often referred to as the {\em Loewner order} \cite{MO} -- is
induced by the cone of all positive semi-definite $n\times n$ matrices: 
$A\le B$ 
(respectively, $A<B$) means that $B-A$ is positive semi-definite (respectively,
positive definite). Using our discussion so far and the arguments further 
below, we can prove the following extension of one of the main results in 
\cite{Lyons}: 

\begin{theorem}\label{th-lyo}
Let $A$ and $B$ be positive contractions such that $A \leq B$ and let 
$\mu$ and $\nu$ be their corresponding determinantal measures. Then 
$\mu \unlhd \nu$ and thus $\mu\pc \nu$. 
\end{theorem}

In \cite[Theorem 8.1]{Lyons} Lyons 
proved the last statement ($\mu\pc \nu$) in Theorem~\ref{th-lyo} under 
the additional hypothesis $[A,B]=0$. 

We will need the following elementary lemma. 

\begin{lemma}\label{l-dado}
If $A,B$ are positive semi-definite $n \times n$ matrices such that 
$A\le B<I$ then $A(I-A)^{-1} \leq B(I-B)^{-1}$.
\end{lemma}

\begin{proof}
By a standard density argument it is enough to prove the lemma for
positive definite matrices. Note first that if $C$ is a positive definite
$n\times n$ matrix such that $C\le I$ then $C^{-1}\ge I$. Moreover, if $D$ is
a positive definite $n\times n$ matrix such that $C\le D$ then 
$M^* CM\le M^* DM$ for any $n\times n$ matrix $M$ since
$$\langle M^* CMu, u\rangle =\langle CMu,Mu\rangle\le \langle DMu,Mu\rangle=
\langle M^* DMu, u\rangle$$
for any $u\in\CC^n$. The latter property with $M=D^{-1/2}$ yields 
$D^{-1/2}CD^{-1/2}\le I$, then by taking inverses and using the first 
property we get $D^{1/2}C^{-1}D^{1/2}\ge I$, and finally by using again the
second property with $M=D^{-1/2}$ we obtain $C^{-1}\ge D^{-1}$. Therefore, 
if $A,B$ are
positive definite matrices as in the lemma it follows
that $I-A^{-1}\le I-B^{-1}$, hence
$$A(I-A)^{-1}=-\left(I-A^{-1}\right)^{-1}\le -\left(I-B^{-1}\right)^{-1}=
B(I-B)^{-1},
$$
which is the desired conclusion.
\end{proof}

\begin{proposition}\label{dado}
Let $A$ and $B$ be positive semi-definite $n \times n$ matrices such that 
$A \leq B$ and set $Z=\diag(z_1,\ldots,z_n)$. Then 
$$
\frac {\det(I+AZ)}{\det(I+A)} \unlhd \frac{\det(I+BZ)}{\det(I+B)} 
\quad \mbox{ and } \quad 
\frac {\det(B+Z)}{\det(B+I)} \unlhd \frac{\det(A+Z)}{\det(A+I)}. 
$$
Moreover, if $A$ and $B$ are positive contractions, then 
$$
\det(I-A+AZ) \unlhd \det(I-B+BZ).
$$
\end{proposition}

\begin{proof}
We just prove the first inequality. The second inequality then follows upon 
considering the operation 
$$\det(I+MZ):=f_M(z_1, \ldots, z_n) \mapsto 
z_1\cdots z_nf_M(z_1^{-1},\ldots, z_n^{-1})$$  
with $M=A$ and $M=B$, respectively, which preserves the set of stable 
polynomials with non-negative coefficients 
(Proposition~\ref{pro-rs}) and obviously reverses the partial order $\unlhd$
(Theorem~\ref{addz}).
Similarly, for the third inequality note that if $A \leq B$ and $A,B$ are 
positive contractions then 
by continuity we may assume that $B<I$. But then  
$A'=A(I-A)^{-1} \leq B(I-B)^{-1}=B'$ by Lemma \ref{l-dado}, hence
$$
\det(I-A+AZ)= \frac {\det(I+A'Z)}{\det(I+A')} \unlhd 
\frac{\det(I+B'Z)}{\det(I+B')}=\det(I-B+BZ). 
$$

We claim that the polynomial  $f \in \RR[z_1,\ldots, z_n,y]$ defined by 
$$f(z,y)=\det(y(I+AZ)+I+BZ):=y^nP_0(z)+ y^{n-1}P_1(z)+\cdots + P_n(z),$$ 
where $z=(z_1,\ldots,z_n)$, is real stable with non-negative coefficients. 
Assuming this claim we see that 
since $P_k(z)= (n-k)!^{-1}\partial^{n-k}f/\partial y^{n-k}\big|_{y=0}$,  
by Remark \ref{dada} we have that 
$$
\det(I+AZ)=P_0 \ll P_1 \ll \cdots \ll P_n=\det(I+BZ), 
$$ 
which then settles the lemma (note that $P_k(0)=\binom n k$, so none of the 
$P_k$'s is identically zero). To prove the remaining claim, note that the 
polynomial 
$
G(z,y)= \det(A+y(B-A)+Z)
$
is stable with non-negative coefficients by Proposition \ref{det-pen}. 
Hence, so are the polynomials 
$$
F(z,y)= z_1\cdots z_ny^nG(z_1^{-1},\ldots, z_n^{-1},y^{-1})
=\det(y(I+AZ)+(B-A)Z) 
$$
and $F(z,y+1)=\det(y(I+AZ)+I+BZ)=f(z,y)$. 
\end{proof}

\begin{proof}[Proof of Theorem~\ref{th-lyo}]
The generating polynomials of $\mu$ and $\nu$ are $\det(I-A+AZ)$ and 
$\det(I-B+BZ)$,  respectively. Combine Propositions~\ref{stdom} and~\ref{dado}.
\end{proof}

We end this section with an open question, namely whether the following
converse to Proposition~\ref{dado} holds. 

\begin{question}\label{q-pos-d}
Let $A$ and $B$ be positive semi-definite $n \times n$ matrices such that 
$$
\frac {\det(I+AZ)}{\det(I+A)} \unlhd \frac{\det(I+BZ)}{\det(I+B)},  
\quad \mbox{ or equivalently, } \quad 
\frac {\det(B+Z)}{\det(B+I)} \unlhd \frac{\det(A+Z)}{\det(A+I)},
$$
where $Z=\diag(z_1,\ldots,z_n)$. Is it true that $A\le B$?
\end{question}

Note that an affirmative answer to Question \ref{q-pos-d} would give a 
characterization of the Loewner order
on the cone of all positive semi-definite $n\times n$ matrices in terms
of determinantal stable polynomials.

\begin{remark}\label{r-spec}
Stochastic domination is closely related to the
notion of majorization (or spectral order) for real vectors 
that was axiomatized by Hardy-Littlewood-P\'olya 
and Schur \cite{HLP,MO}. The latter also induces a partial order (though not 
the same as the proper position) on univariate real-rooted polynomials, 
see, e.g., \cite{Bor}. 
\end{remark}

\subsubsection{Truncations of Strongly Rayleigh Measures and 
Polynomials}\label{sss-trunc}

We will now discuss certain truncations of stable polynomials and strongly 
Rayleigh measures.

\begin{lemma}\label{lum}
Suppose that  $f(z)= \sum_{\alpha\in\NN^n} a(\alpha) z^\alpha \in 
\RR[z_1,\ldots, z_n]$ is stable of total degree at most $d$ and has 
non-negative 
coefficients. For $0\le k\le d$ let 
$$
E_k(z)= \sum_{\mid \alpha \mid = k}a(\alpha)z^\alpha. 
$$
If $0\leq p \leq q \leq d$, then 
\begin{equation}\label{Eqp}
\sum_{k=p}^q \binom {q-p}{k-p} \frac {E_k(z)}{\binom {d}{k}}y^{q-k}
\end{equation}
 is a stable polynomial in the variables  $z_1, \ldots, z_n,y$.
\end{lemma}

\begin{proof}
Define a sequence of polynomials in $\RR[z_1, \ldots, z_n, y]$ as follows: 
\begin{eqnarray*} 
f_1(z,y) &=&  \sum a(\alpha) z^\alpha y^{d-|\alpha|} 
= \sum_{k=0}^d E_k(z)y^{d-k}, \\ 
f_2(z,y) &=& {\partial^ {d-q}f_1}/{\partial y^{d-q}}, \\
f_3(z,y) &=& y^{q}z_1^d\cdots z_n^d 
f_2(z_1^{-1}, \ldots, z_n^{-1}, y^{-1}), \\
f_4(z,y) &=&   {\partial^ {p}f_3}/{\partial y^{p}}, \\ 
f_5(z,y) &=& y^{q-p}z_1^d\cdots z_n^d 
f_4(z_1^{-1}, \ldots, z_n^{-1}, y^{-1}).
\end{eqnarray*}
By Theorem \ref{homhom},  $f_1$ is a stable polynomial, and by the 
closure properties of stable polynomials (Proposition~\ref{pro-rs}) 
we also have that $f_2, \ldots, f_5$ are stable. The polynomial $f_5$ is a 
constant multiple of \eqref{Eqp}. 
\end{proof}

The above lemma gives the following generalization of Newton's inequalities,
which correspond to the special case when $f(z)=(1+z_1)\cdots (1+z_n)$. 

\begin{corollary}\label{Newton-in}
Let $f(z)$ and $E_k(z)$, $k=0,\ldots,d$, be as in the statement of 
Lemma~\ref{lum}. Then 
\begin{equation}\label{N-IE}
\frac {E_k^2(z)}{{\binom {d}{k}}^2} \geq 
\frac {E_{k-1}(z)}{\binom {d}{k-1}}\frac {E_{k+1}(z)}{\binom {d}{k+1}}, 
\quad z \in \RR^n, 1 \leq k \leq d-1.
\end{equation}
\end{corollary}

\begin{proof}
By Lemma~\ref{lum}, the polynomial 
\begin{equation}\label{ert}
 (z,y)\mapsto \frac {E_{k-1}(z)}{\binom {d}{k-1}}y^2+ 
2\frac {E_k(z)}{\binom {d}{k}}y + \frac {E_{k+1}(z)}{\binom {d}{k+1}}
 \end{equation}
 is stable. Hence, for any fixed $z \in \RR^n$ the resulting
(univariate) polynomial in $y$ is real-rooted. The corollary 
 follows upon taking the discriminant in \eqref{ert}. 
 \end{proof}

Recall from Definition~\ref{d-trunc} the truncation $\mu_{p,q}$ 
of a measure $\mu$ to $[p,q]$. It has been discussed whether truncations of 
measures preserve negative dependence properties, cf.~\cite{J-DP,P} and 
references therein. As we show in Corollary \ref{cor-trunc} below, this is 
the case for strongly Rayleigh measures and truncations of length 
$q-p \leq 1$. 
Truncations of length 
$q-p \geq 2$ do not preserve the strong Rayleigh property. However, the 
weighted truncation suggested by \eqref{Eqp} obviously does. 

\begin{corollary}\label{cor-trunc}
Suppose that $\mu$ is a strongly Rayleigh probability measure 
on $2^{[n]}$ and that $0 \leq p \leq q \leq n$ with $q-p \leq 1$. 
Then $\mu_{p,q}$ is strongly Rayleigh. 
\end{corollary}

\begin{proof}
We prove the corollary for $q-p=1$. Let $f$ be the generating polynomial of 
$\mu$ and let $E_k$ be as in the statement of Lemma~\ref{lum}. Then the 
polynomial 
\begin{equation}\label{pq}
 g(z,y)=\frac {E_p(z)}{\binom {n}{p}}y + \frac {E_q(z)}{\binom {n}{q}} 
\end{equation}
is stable. Clearly, the generating polynomial of $\mu_{p,q}$ is a (positive)
constant multiple of $g(z,q/(n-p))$. The corollary follows since stable 
polynomials are closed under setting variables equal to real numbers, 
see Proposition~\ref{pro-rs}. 
\end{proof}

In \S \ref{s-countex} we show that Pemantle's Conjecture~\ref{st} (see 
\S \ref{s24}) fails in its full generality. However, it is true for strongly 
Rayleigh measures: 

\begin{theorem}\label{stst}
Let $\mu$ be a strongly Rayleigh probability measure on $2^{[n]}$, 
and let $1\leq k \leq n$. If 
$\mu(\{S : |S|=k-1\})\mu(\{S:|S|=k\}) \neq 0$, then 
$\mu_{k-1} \pc \mu_k$.  
\end{theorem}

\begin{proof}
Let $f$ be the generating polynomial of $\mu$ and let $E_k$ be as in the 
statement of Lemma~\ref{lum}. By assumption one has 
$E_{k-1}(\ben)E_k(\ben)\neq 0$, where $\ben\in\RR_+^n$ 
is as before the ``all ones 
vector''. From Lemma~\ref{lum} and the last part of Remark~\ref{dada} 
we deduce that the polynomial 
$$
F(z_1,\ldots,z_{n+1})=z_{n+1}
\frac{E_{k-1}(z_1,\ldots, z_n)}{E_{k-1}(\ben)} + 
\frac {E_{k}(z_1,\ldots, z_n)}{E_{k}(\ben)} 
$$
is stable. Since the above quotients are the generating polynomials
of $\mu_{k-1}$ and $\mu_k$, respectively, we have $\mu_{k-1} \unlhd \mu_k$ 
by Theorem \ref{addz} and Definition \ref{d-order}. The desired 
conclusion follows from Proposition \ref{stdom}. 
\end{proof}

\subsection{The Partial Symmetrization Procedure}\label{n-s5}

Recall Remark~\ref{r-complex}, the partial 
symmetrization
of (complex) measures on $2^{[n]}$ defined in \eqref{t-theta}, and
Definition~\ref{d-rs-pol}. The main 
result of this section is the following theorem. 

\begin{theorem}\label{asym}
Let $\mu$ be a complex measure on $2^{[n]}$ with 
(complex) stable generating polynomial. Then for any  
transposition $\tau\in\fS_n$ and $\theta\in[0,1]$ the generating polynomial of
the partial symmetrization 
$\mu^{\tau,\theta}=\theta\mu+(1-\theta)\tau(\mu)$ is (complex) stable.
\end{theorem}

In particular, Theorem~\ref{asym} implies that the class of strongly
Rayleigh measures on $2^{[n]}$ is invariant under the partial symmetrization
procedure. Clearly, the latter result is real stable in nature. However, in 
order
to establish it we have to consider the (wider) complex stable context 
(cf.~Definition~\ref{d-rs-pol}), as in the above theorem. 

For the proof of Theorem~\ref{asym}
we need a multivariate generalization of the classical 
Obreschkoff theorem that was obtained in \cite{BBS2}. 

\begin{theorem}[Multivariate Obreschkoff theorem]
Let $f,g \in \RR[z_1, \ldots z_n]$. Then all non-zero polynomials in the space 
$
\{\alpha f + \beta g: \alpha, \beta \in \RR\}
$
are stable if and only if either $f +ig$ or $f-ig$ is stable or $f=g\equiv 0$. 
\end{theorem}

\begin{remark}
There are now elementary proofs of Theorem \ref{asym} that avoid the 
multivariate Obresckoff theorem, see \cite{BB,L3}. Also, the referee 
sketched an elementary proof similar to those given in \cite{BB,L3}.
\end{remark}



\begin{proof}[Proof of Theorem \ref{asym}] 
Let $f \in \CC[z_1,\ldots,z_n]$ be the generating polynomial of $\mu$. 
Hence $f$ is  stable and multi-affine. 
Assuming, as we may, that $i=1$ and $j=2$, we need to prove that 
$$
\theta f(\xi_1, \xi_2, \ldots, \xi_n)+(1-\theta) 
f(\xi_2, \xi_1,\ldots, \xi_n) \neq 0
$$
whenever $\xi_1, \ldots, \xi_n \in \{z :\Im(z)> 0\}$ 
(cf.~Definition~\ref{d-rs-pol}). By fixing 
$\xi_3,\ldots,\xi_n$ arbitrarily in the open upper half-plane and considering 
the multi-affine polynomial in variables $z_1,z_2$ given by
$$(z_1,z_2)\mapsto g(z_1,z_2):= f(z_1,z_2,\xi_3, \ldots, \xi_n)$$ 
we see that the problem reduces to proving that the linear operator 
$T_\theta$ on multi-affine polynomials in two variables defined by  
$$
T_\theta(h)(z_1,z_2)=\theta h(z_1,z_2)+(1-\theta)h(z_2, z_{1}) 
$$
preserves stability. Let us first show that $T_\theta$ preserves real 
stability. For this let  
$$f(z_1,z_2) = a_{00} + a_{01}z_2 + a_{10}z_1 + a_{11}z_1z_2 \in 
\RR[z_1,z_2]\setminus \{0\}$$ 
be real stable and set $D_f=a_{00}a_{11}-a_{01}a_{10}$. By 
Theorem~\ref{stable-Rayleigh} we know that 
$f$ is real stable if and only if $D_f \leq 0$. Since 
$$
D_{T_\theta(f)} = D_f - \theta(1-\theta)(a_{01}-a_{10})^2\le D_f
$$
it follows that $T_\theta(f)$ is also real stable, as required.
To deal with the (non-real) stable case let first
$h = g +if \in \CC[z_1,z_2]$ be {\em strictly stable} and multi-affine.  
This means that 
$h(z_1,z_2)\neq 0$ whenever $\Im(z_1) \geq 0$ and $\Im(z_2) \geq 0$.  By the 
multivariate Obreschkoff theorem we have that $\alpha g + \beta f$ is real  
stable or zero for all $\alpha, \beta \in \RR$. Hence, by the above, we also 
have that $\alpha T_\theta(f) + \beta T_\theta (g)$ is real stable or zero  
for all $\alpha, \beta \in \RR$. By the  multivariate Obreschkoff theorem 
again this gives that $T_\theta(g+ if)$ or 
$T_\theta(g- if)$ is stable. However, one has 
$$T_\theta(h)(z,z)=h(z,z)=g(z,z)+if(z,z),$$ 
which is strictly stable by assumption. If $T_\theta(g - if)$ is stable then $g(z,z)-if(z,z)$ is stable and 
$g(z,z)+if(z,z)$ is strictly stable which is impossible since it would imply that all zeros of (the strictly stable polynomial)  $g(z,z)+if(z,z)$ were real. Hence   
$T_\theta(g + if)$ is stable. 

Now if $\epsilon>0$ and $p\in\CC[z_1,z_2]$ is stable then clearly 
$p(z_1+i\epsilon,z_2+i\epsilon)$ is strictly stable (in the sense defined 
above) and thus any stable polynomial is the uniform limit on compact sets 
of strictly stable polynomials. By the multivariate version of 
Hurwitz' classical theorem on the ``continuity of zeros'' 
(see, e.g., \cite[Theorem 2.3]{BBS1} and \cite[Footnote 3, p.~96]{COSW}) 
this proves that $T_\theta$ preserves stability. 
\end{proof}

\begin{remark}\label{rem-sym}
From Theorem~\ref{asym} it follows that the symmetrization $\mu_s$ of a 
strongly Rayleigh measure $\mu$ is also strongly Rayleigh. A different 
proof is as follows. If $\mu$ is a strongly Rayleigh measure on $2^{[n]}$ 
then~\eqref{sym} implies that the diagonal
specialization $\De(g_{\mu_s})(t)=\De(g_{\mu})(t)$ (cf.~Definition~\ref{d-ulc})
is a real-rooted univariate polynomial and by the 
Grace-Walsh-Szeg\"o theorem we conclude that the symmetrization
$\mu_s$ of $\mu$ is also strongly Rayleigh.
\end{remark}

\section{Preservation of Negative Dependence by Symmetric Exclusion 
Processes}\label{n-s6}

In this section we study some consequences of the results in \S 
\ref{s-stable} for symmetric exclusion processes and provide an answer to
the problem of finding a natural negative dependence 
property that is preserved by such evolutions (Problem~\ref{p-see}).

Let us first describe these processes in more detail. A configuration of 
particles on a countable set $S$ is a point
$\eta$ on the Boolean lattice $2^S$, which we identify with $\{0,1\}^S$, so 
$\eta(i)=1$ means that site $i$ is occupied, while $\eta(i)=0$ means that it 
is vacant. One is given a non-negative symmetric kernel 
$Q=(q_{i,j})_{i,j\in S}$ satisfying
\begin{equation}\label{eq-1-l}
\sup_i\sum_jq_{i,j}<\infty.
\end{equation}
The corresponding exclusion process $\eta_t$ is the Markov process on 
$\{0,1\}^{S}$ in which
$\eta\rightarrow \tau_{i,j}(\eta)$ at rate $q_{i,j}$, where $\tau_{i,j}$ 
is the transposition
that interchanges the coordinates $\eta(i)$ and $\eta(j)$. Note that if
$\eta(i)=\eta(j)$ (i.e., the two sites are both occupied or both empty), 
this transposition has no effect,
while if $\eta(i)\neq\eta(j)$, the transposition has the effect of moving the
particle from
the occupied site to the vacant site. In the language of semigroups and 
generators, the
semigroup $T(t)$ on $C(\{0,1\}^S)$ defined by 
$T(t)F(\eta)=E^{\eta}F(\eta_t)$ has generator
$$\mathcal{L}F(\eta)=\frac 12\sum_{i,j}q_{i,j}[F(\tau_{i,j}(\eta))-F(\eta)]$$
for functions $F$ that depend on finitely many coordinates. Condition 
\eqref{eq-1-l} guarantees
that the process is well defined and uniquely determined by $Q$.

One can give a more probabilistic description of the symmetric exclusion 
process in
terms of Poisson processes.
Recall that a Poisson process $N(t)$ with rate $\lambda>0$ is a random 
increasing step
function with jumps of size one, in which the times between jumps are 
independent
exponentially
distributed random variables $\tau_k$ with parameter $\lambda$:
$P(\tau_k>t)=e^{-\lambda t}$.
Now take a collection $N_{i,j}(t)$ of independent Poisson processes indexed by
unordered pairs $i,j$ --
$N_{i,j}(t)$ with rate $q_{i,j}$ -- and apply the transposition 
$\tau_{i,j}$ to
$\eta$ at the jump times of
$N_{i,j}(t)$. We will denote probabilities for the process with initial
configuration $\eta$ by $P^{\eta}$.

Negative correlation inequalities have played an important role in the
development of the theory of symmetric exclusion processes. For example,
Proposition 1.7 of \cite[Chap.~VIII]{L0} implies that
\begin{equation}\label{eq-2-l}
P^{\eta}(\eta_t\equiv 1 \text{ on } A)\leq\prod_{i\in
A}P^{\eta}(\eta_t(i)=1),\quad A\subseteq S.
\end{equation}
This inequality was used in an essential way in the characterization of 
stationary distributions
of the process, see Theorem 1.44 in {\em op.~cit} and the end of
this section. Inequality 
\eqref{eq-2-l} was generalized by Andjel in \cite{An} to
\begin{equation}\label{eq-3-l}
\begin{split}
&P^{\eta}(\eta_t\equiv 1\text{ on }A\cup B)\leq P^{\eta}(\eta_t\equiv 1
\text{ on
}A)P^{\eta}(\eta_t\equiv
1\text{ on } B),\\ 
&A,B\subseteq S,\, A\cap B=\emptyset,
\end{split}
\end{equation}
and used to prove an ergodic theorem. For further applications of correlation
inequalities in this
setting, see the references and discussion in \cite{L2}.

In that paper, one of us conjectured that the distribution of $\eta_t$ at 
time $t$ with deterministic initial configuration
$\eta$ is negatively associated. Even the special case
\begin{equation}\label{eq-4-l}
\begin{split}
&P^{\eta}(\eta_t\equiv 1\text{ on }A,\eta_t\equiv 0\text{ on }B)\geq
P^{\eta}(\eta_t\equiv 1\text{ on
}A)P^{\eta}(\eta_t\equiv 0\text{ on } B),\\ 
&A,B\subseteq S,\,A\cap B=\emptyset
\end{split}
\end{equation}
remained open (note that in \eqref{eq-4-l}, one of the events is increasing
while the other is decreasing). 
Surprisingly, despite the similarity between \eqref{eq-3-l} and 
\eqref{eq-4-l}, 
neither the proof of \eqref{eq-3-l} in \cite{An}, nor a somewhat
different one given in \cite{L2}, extends to prove \eqref{eq-4-l}. 
We are now able to prove this conjecture, see Theorem~\ref{th-Lig} below. 
The latter is actually a consequence of a yet stronger property for finite
 symmetric exclusion processes that we will now establish.

\begin{proposition}\label{prop-lig}
Suppose that $S$ is finite and the initial distribution
$\eta_0$ 
of a symmetric exclusion process on $\{0,1\}^S$ is strongly Rayleigh.
Then so is the distribution of $\eta_t$ for all $t> 0$.
\end{proposition}

\begin{proof}  
First, we observe that it suffices to prove the statement in case
$q_{i,j}\neq 0$ for only one
pair $i,j$.  To see this, suppose $T_1(t)$ and $T_2(t)$ are semigroups for 
finite state Markov chains
on the same state space with generators $\mathcal{L}_1$ and $\mathcal{L}_2$ 
respectively, and let $T(t)$ be the semigroup
with generator $\mathcal{L}_1+\mathcal{L}_2$. The Trotter product formula 
(see \cite[p.~33]{EK}) gives
$$T(t)=\lim_{n\rightarrow\infty}\big[T_1(t/n)T_2(t/n)\big]^n.$$
It follows that if a closed set of probability measures $\mathcal{M}$ on 
the state space has the property 
that $\mu\in \mathcal{M}$ implies $\mu T_i(t)\in \mathcal{M}$ for $i=1,2$ 
and $t>0$, then 
$\mu\in\mathcal{M}$ implies $\mu T(t)\in \mathcal{M}$ for $t>0$.
Now write the generator of the symmetric exclusion process as the following 
sum, corresponding to transitions at individual pairs of sites:
$$\mathcal{L}=\sum_{i,j}\mathcal{L}_{i,j},$$
where
$$\mathcal{L}_{i,j}F(\eta)=\frac 12q_{i,j}[F(\tau_{i,j}(\eta))-F(\eta)].$$
Repeated application of the preceding observation with $\mathcal{M}=$ the 
set of probability measures on
$\{0,1\}^S$ with stable generating polynomials, implies that the statement of 
the proposition
is correct for the process with generator $\mathcal{L}$ if it is correct for 
the process with generator $\mathcal{L}_{i,j}$ for each $i,j$. 

In the case of a single non-zero $q_{i,j}$, the generating
polynomial of the distribution of $\eta_t$ is given by $\alpha
f(z)+(1-\alpha)f(\tau_{i,j}(z))$, where
$\alpha=f(N_{i,j}(t)\text{ is even})$ and $f$ is the generating polynomial
for the initial distribution. Therefore, the result follows from 
Theorem~\ref{asym}.
\end{proof}

Recall from \S \ref{s240} that if $S$ is a 
countably infinite set, a probability measure on
$\{0,1\}^S$ is strongly
Rayleigh if every projection onto a finite subset of $S$ is strongly Rayleigh. 
Note that this property is preserved by weak convergence of measures. 

\begin{theorem}\label{th-Lig}
Suppose that $S$ is countable and the initial distribution $\eta_0$ of a 
symmetric exclusion process 
on $\{0,1\}^S$ is strongly Rayleigh. Then the distribution of $\eta_t$ is 
strongly Rayleigh, and hence CNA+, for all $t>0$.
\end{theorem}

\begin{proof} 
Since stability (the strongly Rayleigh property) implies
negative association by Theorem~\ref{stab-NA}, the result follows from 
Proposition~\ref{prop-lig} when $S$ is finite. 
For general $S$, approximate the exclusion process on $S$ by processes on 
finite subsets $S_n$ of
$S$ that increase to $S$, by taking the transposition rates for the process 
on $S_n$ to be $q_{i,j}$ if
$i,j\in S_n$ and zero otherwise. By the Trotter-Kurtz semigroup convergence 
theorem (see \cite[p.~28]{EK}), the
processes on $S_n$ converge to the process on $S$, so the result holds on 
$S$ as well. (See \cite[\S I.3]{L0} for details.) 
\end{proof}

\begin{remark}\label{prod-meas}
In particular, Theorem \ref{th-Lig} proves Conjecture \ref{con-lig} and shows 
that its conclusion remains valid if the initial configuration is a 
product measure (since any such measure is obviously strongly Rayleigh, 
cf.~Definition \ref{d-prod-m} and \S \ref{s240}).
\end{remark}

We conclude this section with an application of Theorem~\ref{th-Lig} 
to the extremal stationary
distributions for an irreducible
symmetric exclusion process. Theorem~1.44 of \cite[Chap.~VIII]{L0} states
that these
are exactly $\{\mu_{\alpha}:\alpha\in \mathcal{H}\}$, where 
\begin{itemize}
\item[(a)] $\mathcal{H}=\left\{\alpha:S\rightarrow [0,1] : 
\sum_jq_{i,j}[\alpha(j)-\alpha(i)]=0\text{ for each }i\right\}$ are the
harmonic functions for $Q$ with values in $[0,1]$, and 
\item[(b)] $\mu_{\alpha}=\lim_{t\rightarrow\infty}\nu_{\alpha}T(t)$, 
where $\nu_{\alpha}$ is the product measure on $\{0,1\}^S$ satisfying 
$$\nu_{\alpha}\{\eta:\eta\equiv 1\text{ on }A\}=\prod_{i\in A}\alpha(i),\quad
A\subseteq S.$$
\end{itemize} 
Very little is known about $\mu_{\alpha}$ other than 
\begin{itemize}
\item[(I)] 
$\mu_{\alpha}=\nu_{\alpha}$ if and only if $\alpha$ is constant, and 
\item[(II)] 
$\mu_{\alpha}\{\eta:\eta(i)=1\}=\alpha(i)$ for each $i\in S$. 
\end{itemize}
As a consequence of Theorem~\ref{th-Lig}, we can now add
\begin{itemize}
\item[(III)] 
$\mu_{\alpha}$ is negatively associated for each $\alpha\in\mathcal{H}$. 
\end{itemize}

\begin{remark}\label{r-lig-21}
If $q_{i,j}>0$ for all $i,j$, then as we noted in \S \ref{ss-ex-case} 
the limiting distribution of $\eta_t$ as
$t\to\infty$ is the symmetrization of the initial distribution. In 
\S \ref{s-countex} we construct examples of measures in $\fP_{n}$, 
$n\ge 20$, that are CNA+ but do
not have ULC rank sequences. By Theorem~\ref{pro-w} in \S \ref{ss-ex-case}, 
these counterexamples
to Conjectures~\ref{c-p}--\ref{cna-ulc} show that none of 
the properties NLC, h-NLC, Rayleigh/h-NLC+, CNA, CNA+ is preserved by the
symmetrization procedure \eqref{sym}, and by \cite[Theorem 3.2]{L2}, the same
can be said about the NA property. Therefore, the analog of 
Proposition~\ref{prop-lig} fails for all the aforementioned properties.
\end{remark}

\begin{remark}\label{r-lig-2}
The analog of Theorem~\ref{th-Lig} for asymmetric exclusion processes is 
false. A simple example is
obtained by taking $S=\{1,2\}$ in which only transitions from state 1 to 
state  2 are allowed. If the initial
distribution gives probability $\frac 14$ to each of the configurations in
$\{0,1\}^S$, then
the limiting distribution as $t\to\infty$ is given by 
$\mu(\{11\})=\mu(\{00\})=\frac14$,
$\mu(\{10\})=0$, and $\mu(\{01\})=\frac 12$. This measure is not NA.
\end{remark}

\section{Almost Exchangeable Measures}\label{Almost}

By Pemantle's result (Theorem~\ref{pro-w} in \S \ref{ss-ex-case}), for 
symmetric measures in $\fP_n$ all five conditions NLC, h-NLC, 
Rayleigh/h-NLC+, CNA, CNA+ are equivalent to the ULC condition for the rank 
sequence. Therefore, the relations between 
all these negative dependence properties are completely understood in the 
exchangeable case. In this section we study the almost exchangeable
case, that is, measures whose generating polynomials are 
symmetric in all but possibly one variable (Definition~\ref{d-exc}). We prove 
that for such measures one has both 
h-NLC $\Leftrightarrow$ CNA and 
Rayleigh/h-NLC+ $\Leftrightarrow$ CNA+, which confirms 
Pemantle's Conjecture \ref{con-cna-R} in this case. In the next section we 
will use this result to show that the counterexamples in 20 or more 
variables  that we 
construct there satisfy 
all negative dependence properties other than strongly Rayleigh.

Let $C$ be an $m \times n$ matrix (where $m,n$ are allowed to be infinite), 
and let $\mathbf{i}$ and $\mathbf{j}$ be 
finite equi-numerous subsets of $[m]$ and $[n]$, respectively. Denote by 
$C(\mathbf{i}, \mathbf{j})$ the minor of $C$ with rows indexed by 
$\mathbf{i}$ and columns indexed by $\mathbf{j}$. The Cauchy-Binet theorem, 
see, e.g., \cite{Kar}, asserts that 
$$
(AB)(\mathbf{i}, \mathbf{j})
= \sum_{\mathbf{k},\,|\mathbf{k}|=|\mathbf{i}|=|\mathbf{j}|} A(\mathbf{i}, 
\mathbf{k})B(\mathbf{k}, \mathbf{j})
$$
for any $m\times n$ matrix $A$ and $n\times q$ matrix $B$, where
$1\le m,n,q\le \infty$, provided that the sum in the right-hand side converges
absolutely. 

Recall the notion of (non-negative) $\LC$ sequence from 
Definition \ref{d-ulc}. For convenience, we first prove a technical result 
that will be frequently used later on.

\begin{lemma}\label{l-tech}
Let $\{a_j\}_{j=0}^n$ and $\{b_j\}_{j=0}^n$ be two LC 
sequences that satisfy
\begin{itemize} 
\item[(i)] $\{a_j +  b_j\}_{j=0}^n$ is LC, and 
\item[(ii)] $a_jb_{j+1} \geq a_{j+1}b_j$ for all $0 \leq j \leq n-1$. 
\end{itemize}
If $\ell<k$ and $b_\ell a_k>0$, then $b_ra_r>0$ for all $\ell\leq r\leq k$.
\end{lemma}

\begin{proof}
Suppose that $a_r=0$ for
some $\ell \leq r< k$ and let $s$ the largest such index. Assumption (ii) 
for $j=s$ gives $b_s=0$, hence $a_s+b_s=0$. This contradicts the 
``no internal zeros" part of assumption (i) (cf.~Definition~\ref{d-ulc}). 
Therefore $a_r>0$, and a similar argument shows that $b_r>0$ as well. 
\end{proof}

The following proposition shows that for almost exchangeable measures the 
Rayleigh property is the same as the negative lattice condition (NLC) 
closed under conditioning and external fields on the last variable. 
Indeed, condition (ii) follows directly from 
NLC while (i) follows from NLC after having imposed an external field on the 
last coordinate and then projected onto the first $n$ coordinates. 

\begin{proposition}\label{RAS}
Let $\mu\in\fP_{n+1}$ be such that its generating polynomial $g=g_\mu$ is 
symmetric in its first $n$ variables, so that
$$
g(z_1,\ldots,z_{n+1})= z_{n+1}\sum_{k=0}^{n}a_ke_k(z_1,\ldots, z_n)+
\sum_{k=0}^{n}b_ke_k(z_1,\ldots, z_n).
$$
Then $\mu$ is Rayleigh if and only if 
\begin{itemize} 
\item[(i)] $\{\theta a_k + (1-\theta) b_k\}_{k=0}^n$ is LC 
for all $0 \leq \theta \leq 1$, and 
\item[(ii)] $a_kb_{k+1} \geq a_{k+1}b_k$ for all $0 \leq k \leq n-1$. 
\end{itemize}
\end{proposition}
\begin{proof}
Suppose first that $\mu$ is Rayleigh. By Proposition~\ref{p-proper} we get 
that if
$0 \leq \theta < 1$ then
$g(z_1, \ldots, z_n,\theta/(1-\theta))$ is a Rayleigh polynomial that is 
symmetric in all its $n$ variables. Condition (i) now follows from the 
exchangeable case (Theorem~\ref{pro-w}).  
Let $S \subseteq [n-1]$ be of cardinality $i$.
The NLC condition gives 
$$
a_ib_{i+1} = \mu(S \cup \{n+1\})\mu(S \cup \{n\}) \geq 
\mu(S \cup \{n,n+1\})\mu(S)= a_{i+1}b_i, 
$$
which verifies (ii). 

Assume now that $\{a_k\}_{k=0}^{n}$ and $\{b_k\}_{k=0}^n$ satisfy (i) and 
(ii). The inequality 
$$
\frac {\partial g} {\partial z_i}(x) \frac {\partial g} {\partial z_j}(x) 
\geq \frac {\partial^2 g} {\partial z_i \partial z_j}(x) g(x), 
\quad x\in \RR_+^{n+1}, \,i,j \in [n],
$$
follows from the exchangeable case and (i). What is left to prove is the 
inequality 
$$
\frac {\partial g} {\partial z_i}(x) \frac {\partial g} {\partial z_{n+1}}(x) 
\geq \frac {\partial^2 g} {\partial z_i \partial z_{n+1}}(x) g(x), 
\quad x\in \RR_+^{n+1},\, i\in [n], 
$$
and by symmetry (in the first $n$ variables) we only need to check it for
e.g.~$i=n$. We may write $g$ as 
\begin{multline*}
g(z_1,\ldots,z_{n+1})= \sum_{k=0}^{n} b_k e_k(z_1, \ldots, z_{n-1})+
z_n\sum_{k=0}^{n} b_{k+1} e_k(z_1, \ldots, z_{n-1}) \\
+ z_{n+1}\sum_{k=0}^n a_k e_k(z_1, \ldots, z_{n-1})+ z_nz_{n+1} 
\sum_{k=0}^n a_{k+1} e_k(z_1, \ldots, z_{n-1}),
\end{multline*}
so the desired inequality is equivalent to 
\begin{equation}\label{eqq}
\left(\sum_j a_je_j\right)\!\left(\sum_j b_{j+1}e_j\right)
-\left(\sum_j a_{j+1}e_j\right)\!\left(\sum_j b_je_j\right)\geq 0
\end{equation}
for all $x \in \RR_+^{n-1}$, where $e_j=e_j(x_1,\ldots, x_{n-1})$. Let now
 $$C=\left( \begin{array}{ccccc}
  a_0 & a_1 & a_2 & \ldots & a_n \\
  b_0 & b_1 & b_2 & \ldots  & b_n
\end{array} \right)$$
and recall the definition of a $\TP_p$ matrix from \S \ref{ss-ex-case}.
We claim that $C$ is $\TP_2$. Let $1 \leq i<j \leq n$. We need to 
prove that 
$a_ib_j \geq b_ia_j$. Clearly, we may assume that $b_ia_j \neq 0$ and that 
$j-i \geq 2$. By Lemma~\ref{l-tech} we have $b_ra_r>0$ for $i\le r\le j$ and
thus 
$$
\frac {a_i} {b_i} \geq \frac {a_{i+1}} {b_{i+1}} \geq \cdots 
\geq \frac {a_j} {b_j}
$$
by (ii), which shows that $C$ is indeed $\TP_2$. 
Let $E=(e_{i-j})_{i,j=0}^{n}$. By Theorem~\ref{ESC} we have that $E$ is 
$\TP_2$ and by the Cauchy-Binet formula so is $CE$.  Now 
 $$CE=\left( \begin{array}{ccc}
  \sum a_je_j  & \sum a_{j+1}e_j\ & \ldots  \\
  \sum b_je_j & \sum b_{j+1}e_j & \ldots 
\end{array} \right),$$
from which \eqref{eqq} follows. 
\end{proof}

The next result is the analog of Proposition~\ref{RAS} for the h-NLC property.

\begin{proposition}\label{pro-h-NLC}
Let $\mu\in \fP_{n+1}$ be such that its  generating polynomial $g=g_{\mu}$
is symmetric in its first $n$ variables, so that
$$g(z_1,\ldots,z_{n+1})
=z_{n+1}\sum_{k=0}^n a_ke_k(z_1,\ldots,z_n)+
\sum_{k=0}^n b_ke_k(z_1,\ldots,z_n).$$
Then $\mu$ is $h$-NLC if and only if
\begin{itemize}
\item[(i)] $\{a_k\}_{k=0}^n$, $\{b_k\}_{k=0}^n$, and 
$\{a_k+b_k\}_{k=0}^n$ are LC, and
\item[(ii)] $a_kb_{k+1}\geq a_{k+1}b_k$ for all $0\leq k\leq n-1$. 
\end{itemize}
\end{proposition}

\begin{proof}
Let $S,T\subset [n]$. The NLC condition applied to each of the three pairs
$(S,T)$, $(S,T\cup\{n+1\})$, and $(S\cup\{n+1\},T\cup\{n+1\})$ gives  
respectively
$$b_{|S|}b_{|T|}\geq b_{|S\cup T|}b_{|S\cap T|},\quad b_{|S|}a_{|T|} 
\geq a_{|S\cup T|}b_{|S\cap T|},
\quad a_{|S|}a_{|T|}\geq a_{|S\cup T|}a_{|S\cap T|}.$$
Therefore, setting $\ell=|S|$, $k=|T|$, $m=|S\cap T|$, we see that 
NLC is  equivalent to
\begin{equation}\label{tech1}
b_\ell b_k\geq b_{k+\ell-m}b_m,\quad b_\ell a_k\geq a_{k+\ell-m}b_m,\quad 
a_\ell a_k \geq a_{k+\ell-m}a_m
\end{equation}
for $m\leq k,\ell\leq n$ and $k+\ell\leq n+m$.

If $S\subseteq [n]$ has cardinality $n-j$, the measure obtained by  
projecting $\mu$ onto $2^{S\cup\{n+1\}}$
is of the same form as $\mu$ with  coefficients
$$a_k(j)=\sum_{i=0}^j\binom jia_{k+i}\quad\text{and}\quad 
b_k(j)=\sum_ {i=0}^j\binom jib_{k+i},
\quad 0\leq k\leq n-j,$$
while the one obtained by projecting $\mu$ onto $2^S$ is exchangeable with  
coefficients $a_k(j)+b_k(j)$, $0\leq k\leq n-j$.
Therefore, h-NLC is equivalent to the statement that for each 
$0 \leq j\leq n$ one has 
\begin{equation}\label{tech2}
\begin{split}
&b_\ell(j)b_k(j)\geq b_{k+\ell-m}(j)b_m(j),\\  
&b_\ell(j)a_k(j)\geq  a_{k+\ell-m}(j)b_m(j),\,\,
a_\ell(j)a_k(j)\geq a_{k+\ell-m}(j)a_m(j),\\
&(a_\ell(j)+b_\ell(j))(a_k(j)+b_k(j))\geq 
(a_{k+\ell-m}(j)+b_{k+\ell-m}(j))(a_m(j) +b_m(j))
\end{split}
\end{equation}
for $m\leq k,\ell\leq n-j$ and $k+\ell\leq n-j+m$.

It follows from this that the h-NLC property implies that the sequences in 
(i) have no internal zeros.
For example, if $k\leq \ell$, $b_{k-1}>0, b_k=\cdots=b_\ell=0$, 
$b_{\ell+1}>0$  contradicts $b_kb_\ell\geq
b_{k-1}b_{\ell+1}$. Given this, the fact that the three sequences are LC  
and satisfy (ii) is
a consequence of the first and last inequalities in \eqref{tech1} with 
$\ell=k$, $m=k-1$, the middle
inequality in \eqref{tech1} with $\ell=k+1$, $m=k$, and the last inequality 
in \eqref{tech2} with $j=0$, $\ell=k$, $m=k-1.$

For the converse, assume that (i) and (ii) hold. Note that for  
each $j$, $\{a_k(j)\}_{k=0}^{n-j}$
is a convolution of two LC sequences, and hence is LC by 
\cite[Theorem 1.2 on p.~394]{Kar},
and the same holds for the sequences $b_k(j)$ and $a_k(j)+b_k(j)$.  
This gives three of the
four inequalities in \eqref{tech2}. For the other one, note that since  
evaluating the elementary symmetric
polynomials at the ``all ones vector'' produces binomial coefficients, 
it follows 
from \eqref{eqq} that
\begin{equation}\label{tech3}
a_k(j)b_{k+1}(j)\geq a_{k+1}(j)b_k(j).
\end{equation}
We need to check that 
$$a_k(j)b_\ell(j)\geq a_{k+\ell-m}(j)b_m(j)$$
for $m\leq k,\ell$, and in doing so, can assume that the right-hand side is  
not zero. By Lemma~\ref{l-tech},
it follows that $a_i(j)b_i(j)>0$ for $m\leq i\leq k+\ell-m$. To conclude  
the proof, write
$$\frac{a_k(j)}{a_{k+\ell-m}(j)}\geq 
\frac{a_m(j)}{a_\ell(j)}\geq\frac{b_m (j)}{b_\ell(j)}.$$
The first inequality comes from the already proved part of 
\eqref{tech2}, while  the second
is obtained by repeated application of \eqref{tech3}. 
\end{proof}

\begin{theorem}\label{th-a-ex}
Let $\mu\in\fP_{n+1}$ be such that its generating polynomial $g=g_\mu$ 
is symmetric in its first $n$ variables, so that
$$
g(z_1,\ldots,z_n)= z_{n+1}\sum_{k=0}^{n}a_ke_k(z_1,\ldots, z_n)+
\sum_{k=0}^{n}b_ke_k(z_1,\ldots, z_n).
$$
Suppose that $\{a_k\}_{k=0}^n$ and $\{b_k\}_{k=0}^n$ are LC 
sequences that satisfy
\begin{itemize} 
\item[(i)] $\{a_k +  b_k\}_{k=0}^n$ is LC, and 
\item[(ii)] $a_kb_{k+1} \geq a_{k+1}b_k$ for all $0 \leq k \leq n-1$. 
\end{itemize}
Then $\mu$ is $\CNA$. 
\end{theorem}

\begin{proof}
As in \S \ref{s21}, we denote by $X_i$, $i\in [n+1]$, the $i$-th coordinate
function/random (binary) variable on $2^{[n+1]}$. First, observe that it is 
enough to prove that $\mu$ is  NA. To see this, note
that if $\mu$ is conditioned on 
$X_i=1$ for $k$ values of $i\leq n$, on $X_i=0$ for
$l$ values of $i\leq n$, and possibly on the value of $X_{n+1}$, the  
resulting measure is
of the same form with a new value $n'$ of $n$ and new $a_i$'s and $b_i$'s  
that again satisfy the assumptions of the
theorem. Using primes to denote the new values, we have $n'=n-k-l$ and if the  
value of $X_{n+1}$ is not conditioned on,
then, except for a constant factor, $a_i'=a_{i+k}$ and $b_i'=b_{i+k}$. If
the value of $X_{n+1}$ is conditioned on, then $a_i'=b_i'=a_{i+k}$ if  
the conditioning is on
$X_{n+1}=1$ and $a_i'=b_i'=b_{i+k}$ if the conditioning is on $X_{n+1} =0$.
Now take two increasing functions $F$ and $G$ on $2^{[n+1]}$ that  depend on 
disjoint sets of coordinates (cf.~Definition~\ref{d-inc-ev}). Without loss of 
generality, we may assume that there  is an $m$ satisfying
$1\leq m\leq n$ such that $F$ depends on the coordinates $\{X_1\ldots,X_m \}$, 
and $G$ depends on the
coordinates $\{X_{m+1},\ldots,X_{n+1}\}$. We need to prove
\begin{equation}\label{cna1}
\int FGd\mu\leq\int Fd\mu\int G d\mu.
\end{equation}
Since $\mu$ is symmetric in the first $n$ coordinates, these three  
integrals are not
changed if $F$ and $G$ are replaced by the functions obtained by  
symmetrizing them with
respect to the coordinates $\{X_1,\ldots,X_m\}$ and $\{X_{m+1},...,X_n\}$, 
respectively. These symmetrized
functions are also increasing. Therefore, we may assume that $F$ and  $G$ are 
of the following form:
\begin{equation*}
\begin{split}
&F(X_1,\ldots,X_m)=f_k\text{ if }\sum_{i=1}^mX_i =k,\\
&G(X_{m+1},\ldots,X_n,X_{n+1}=0)=g_k\text{ if }\sum_{i=m+1}^nX_i =k,
\text{ and}\\
&G(X_{m+1},\ldots,X_n,X_{n+1}=1)=h_k\text{ if }\sum_{i=m+1}^nX_i =k,
\end{split}
\end{equation*}
where $f_k,g_k$ and $h_k$ are increasing, and $g_k\leq h_k$ for each  $k$.
Using the fact that the covariance of $F$ and $G$ with respect to the  
measure $\mu$
can be written as
$$\frac 12\sum_{S,T\in 2^{[n+1]}}[F(S)-F(T)][G(S)-G(T)]\mu(S)\mu(T),$$
\eqref{cna1} becomes
\begin{multline*}
\sum_{i,j,k,\ell}\binom mi\binom mj\binom {n-m}k\binom {n-m}\ell 
(f_j-f_i)\bigg[(h_k-h_\ell)a_{k+i}a_{\ell+j}\\+
(g_k-g_\ell)b_{k+i}b_{\ell+j}+(h_k-g_\ell)a_{k+i}b_{\ell+j}
+(g_k-h_\ell)a_{\ell+j}b_{k +i}\bigg]\geq 0.
\end{multline*}
Note that we are using the usual convention that $\binom nm=0$ unless  
$0\leq m\leq n$.
Since the summand above is not changed if the roles of $i$ and $j$  
are interchanged and the roles
of $k$ and $\ell$ are interchanged, this can be written as
\begin{multline}\label{cna2}
\sum_{i<j}\sum_{k,\ell}\binom mi\binom mj\binom {n-m}k\binom
{n-m}\ell(f_j-f_i)\bigg[(h_k-h_\ell)a_{k+i}a_{\ell+j}\\
+ (g_k-g_\ell)b_{k+i}b_{\ell +j}+(h_k-g_\ell)a_{k+i}b_{\ell+j}
+(g_k-h_\ell)a_{\ell+j}b_{k+i}\bigg]\geq 0.
\end{multline}
So, it is enough to show that for fixed $i<j$, and fixed $p$,
\begin{multline*}
\sum_{k+\ell=p}\binom {n-m}k\binom {n-m}\ell
\bigg[(h_k-h_\ell)a_{k+i} a_{\ell+j}\\
+ (g_k-g_\ell)b_{k+i}b_{\ell+j}+(h_k-g_\ell)a_{k+i}b_{\ell+j}
+(g_k-h_\ell)a_{\ell+j}b_{k +i}\bigg]\geq 0.
\end{multline*}
Furthermore, it is enough to prove this in case when the $g_k$'s and $h_k $'s 
take only the values 0 and 1,
since any increasing function on a partially ordered set (in this  case, 
$\{0,\ldots,n\}\times \{0,1\}$)
can be written as a positive linear combination of increasing  
functions that take only the values 0 and 1.

We will in fact prove that for such $g_k$'s and $h_k$'s,
\begin{multline*}
S_q=\sum_{k+\ell=p,\,|k-\ell|\leq 2q}
\binom{n-m}k\binom{n-m}\ell\bigg[(h_k-h_\ell)a_{k+i}a_{\ell+j}\\
+ (g_k-g_\ell)b_{k+i}b_{\ell+j}+(h_k-g_\ell)a_{k+i}b_{\ell+j}
+(g_k-h_\ell)a_{\ell+j}b_ {k+i}\bigg]\geq 0
\end{multline*}
by induction on $q$.  For the
basis step, note that $S_0=0$ if $p$ is odd, while if $p=2r$ is even,  
then
$$S_0=\binom{n-m}{r}^2(h_r-g_r)(a_{r+i}b_{r+j}-a_{r+j}b_{r+i}),$$
which is non-negative by (ii), since $i<j$.

We will prove now that $S_q\geq 0$, assuming that $S_{q-1}\geq 0$.  
The difference $S_q-S_{q-1}$ consists
of two summands, which when combined, become
\begin{multline}\label{cna3}
\binom {n-m}k\binom {n-m}\ell\bigg[(h_k-h_\ell)(a_{k+i}a_{\ell+j}
-a_ {\ell+i}a_{k+j})+
(g_k-g_\ell)(b_{k+i}b_{\ell+j}-b_{\ell+i}b_{k+j})\\
+(h_k-g_\ell)(a_{k+i}b_{\ell+j}
-a_ {k+j}b_{\ell+i})+
(g_k-h_\ell)(a_{\ell+j}b_{k+i}-a_{\ell+i}b_{k+j})\bigg],
\end{multline}
where $k+\ell=p$ while $k-\ell=2q$ if $p$ is even and 
$k-\ell=2q-1$ if $p$ is odd.
Now we consider the possible values of $g_\ell,h_k$ for this pair $k,\ell$.  
Recall that
$g_\ell\leq h_\ell,g_k\leq h_k$. If $g_\ell=h_k$, then $g_t=h_t=$ this common  
value for all $\ell\leq t\leq k$,
so that \eqref{cna3} is zero. So, we may assume that $g_\ell=0,h_k=1$. 
Now there  are four
possible values for the pair $(h_\ell,g_k)$ (note that the use of the 
fractions that appear below is justified by Lemma~\ref{l-tech}):

\medskip

(a) $\quad h_\ell=g_k=1$. In this case, \eqref{cna3} becomes
$$\binom {n-m}k\binom
{n-m}\ell\bigg[(b_{k+i}b_{\ell+j}-b_{\ell+i}b_{k+j})
+(a_{k+i}b_{\ell+j}-a_{k+j}b_ {\ell+i})\bigg].$$
This is non-negative by the log-concavity of the $b$'s and (ii):
$$\frac{a_{k+i}}{a_{k+j}}\geq\frac{b_{k+i}}{b_{k+j}}\geq
\frac{b_{\ell+i}} {b_{\ell+j}}.$$

\medskip

(b) $\quad h_\ell=0,g_k=1$. Now \eqref{cna3} is
$$\binom {n-m}k\binom {n-m}\ell\bigg[(a_{k+i}+b_{k+i})(a_{\ell+j}+b_{\ell+j})
- (a_{\ell+i}+b_{\ell+i})(a_{k+j}+b_{k+j})
\bigg],$$
which is non-negative by (i).

\medskip

(c) $\quad h_\ell=g_k=0$. This time \eqref{cna3} is
$$\binom {n-m}k\binom {n-m}\ell\bigg[(a_{k+i}a_{\ell+j}-a_{\ell+i}a_{k+j})+
(a_{k+i}b_{\ell+j}-a_{k+j}b_{\ell+i})\bigg],$$
which is non-negative by the log-concavity of the $a$'s and (ii).

In these three cases, \eqref{cna3} is non-negative, and therefore 
$S_q \geq 0$ by the induction
hypothesis. The fourth case is different, since \eqref{cna3} need not be  
non-negative:

\medskip

(d) $\quad h_\ell=1, g_k=0$. In this case, $g_t\equiv 0$ and $h_t\equiv 1$ 
for  
$\ell\leq t\leq k$, so
$$S_q=\sum_{r+s=p,\,|r-s|\leq 2q}\binom {n-m}r\binom  {n-m}s(
a_{r+i}b_{s+j}-a_{s+j}b_{r+i}).$$
Now use summation by parts, together with the fact that
$$u(r)=\binom{n-m}r\binom{n-m}{p-r}$$
is decreasing in $r$ for $r\geq p/2$, to argue that it suffices to  show that
\begin{equation}\label{cna4}
T_q=\sum_{r+s=p,\,r,s\geq 0,\,|r-s|\leq 2q}(
a_{r+i}b_{s+j}-a_{s+j}b_{r+i})\geq 0
\end{equation}
for each $q$. Here is the argument: noting that
$$S_q-S_{q-1}=(T_q-T_{q-1})u\big(q+[p/2]\big)$$
and summing, we have
\begin{equation*}
\begin{split}
S_q=&\,\,\sum_{r\leq q}(S_r-S_{r-1})=\sum_{r\leq q}(T_r-T_{r-1})u 
\big(r+[p/2]\big)\\
=&\,\,T_qu\big(q+[p/2]\big)+\sum_{r<q}T_r\big[u\big(r+[p/2]\big)
-\big[u\big(r+1+[p/2]\big)\big].
\end{split}
\end{equation*}
Now to check \eqref{cna4}, we rearrange the sum as follows, 
noting  that there
may be some cancellation of terms, leading to the truncation of the  
upper limit
of the sum:
$$\sum_{r=\ell+i}^{(k+i)\wedge(\ell+j-1)}(a_rb_{p-r+i+j}-a_{p-r+i+j}b_r).$$
Here again $k+\ell=p,$ and $k-\ell=2q$ if $p$ is even and $k-\ell=2q-1$ if 
$p$  is odd. The summands in this sum are
non-negative by (ii) since for all $r$ appearing in the sum, $r\leq p-r +i+j$.
To check this, write $2r=r+r\leq (k+i)+(\ell+j)=p+i+j$.
\end{proof}

\begin{corollary}\label{cor-cna+}
Suppose $\{a_k\}_{k=0}^n$ and $\{b_k\}_{k=0}^n$ are LC sequences
that satisfy
\begin{itemize}
\item[(i)] $\{\lambda a_k+b_k\}_{k=0}^n$ is LC for all $\lambda>0$, and
\item[(ii)] $a_kb_{k+1}\geq a_{k+1}b_k$ for $0\leq k<n.$
\end{itemize}
Then $\mu$ is CNA+.
\end{corollary}

\begin{proof}
Let $\{w_i>0\}_{i=1}^{n+1}$ be the external field that is  applied to $\mu$, 
and
$\mu'$ be the resulting measure:
$$\mu'(S)=\text{ constant }\times\prod_{i=1}^{n+1}w_i^{X_i(S)}\mu(S), \quad 
S\subseteq [n+1].$$
We need to show that \eqref{cna1} holds with $\mu$ replaced by $\mu'$. 
First,  note that the effect
of $w_{n+1}$ is to replace $a_k$ by a constant multiple of $a_k$, so  
since we have introduced a general
positive $\lambda$ in assumption (i), we may assume that $w_{n+1}=1.$  
The analog of \eqref{cna2} that we must check is
\begin{multline}\label{cna5}
\sum_{i<j}\sum_{k,\ell}(f_j-f_i)c_ic_jd_kd_\ell
\bigg[(h_k-h_\ell)a_ {k+i}a_{\ell+j}\\+
(g_k-g_\ell)b_{k+i}b_{\ell+j}+(h_k-g_\ell)a_{k+i}b_{\ell+j} 
+(g_k-h_\ell)a_{\ell+j}b_{k +i}\bigg]\geq 0,
\end{multline}
where now the $f_k$'s, $g_k$'s and $h_k$'s are obtained by  symmetrizing $F$ 
and $G$ with respect to the
field on the set of coordinates $\{X_1,\ldots,X_m\}$ in the case of $F$  
and $\{X_{m+1},\ldots,X_n\}$
in the case of $G$. For example,
$$f_k=c_k^{-1}\cdot\sum F(S)\prod_{i=1}^mw_i^{X_i(S)},$$
where 
$$c_k=\sum \prod_{i=1}^mw_i^{X_i(S)}$$
and the above sums are taken over all $S\subseteq [m]$ such that 
$\sum_{i=1} ^mX_i(S)=k$. The formulas for $g_k$ and $h_k$ are exactly 
analogous and the corresponding normalization factor (instead of $c_k$) is
$$d_k=\sum \prod_{i=m+1}^{n}w_{i}^{X_i(S)},$$
the sum being taken over all $S\subseteq \{m+1,\ldots,n\}$ such that 
$\sum_{i=m+1}^{n}X_i(S)=k$. 
Note that if $w_i\equiv 1$, then $c_k=\binom mk$ and $d_k=\binom{n-m} {k}$, 
so that $f_k$ reduces to the symmetrization used in the proof of 
Theorem~\ref{th-a-ex}, and \eqref{cna5} reduces
to \eqref{cna2}.

In order to deduce \eqref{cna5} from the arguments used to prove \eqref{cna2},
we then  need to check the following statements:

\medskip

(a) $\quad f_k\le f_{k+1}$, $g_k\le g_{k+1}$, $h_k\le h_{k+1}$ and
$g_k\leq h_k$ for all $k$, and

(b) $\quad \{d_k\}_{k=0}^{n-m}$ is an LC sequence.

\medskip

\noindent The latter is needed to check that
$d_rd_{p-r}$
is decreasing in $r$ for $r\geq p/2$, which is then used in the analog  
of the treatment of case (d) in the proof of Theorem~\ref{th-a-ex}. Note that 
statement (b) follows from e.g.~Theorem~\ref{ESC}, since $d_k$ is just the
$k$-th elementary symmetric function on $\{w_{m+1},\ldots,w_n\}$.
 
Now the fact that $g_k\leq h_k$ is immediate from the monotonicity of $G$  
in its last argument. The
verification of the first three statements in (a) is the  same 
for all of them, so
we check (a) only for $f_k$. Letting as usual $|S|$ denote the cardinality
of a set $S$, write
$$c_kc_{k+1}(f_{k+1}-f_k)= 
\sum_{|S|=k,\,|T|=k+1} [F(T)-F(S)] \prod_{i=1}^mw_i^{X_i(S)+X_i(T)}.$$
Consider the sum of the terms corresponding to $S$'s and $T$'s that  
satisfy $X_i(S)+X_i(T)=z_i$
for all $i$ and an arbitrarily fixed $z\in\{0,1,2\}^m$. If $A,B,C$ are the  
index sets corresponding to those coordinates $X_i$ for which
$z_i=0,1,2$, respectively, we see that we need to check
\begin{equation}\label{cna6}
\sum_{|T|=k+1,\,X_i(T)=0,\,i\in A,\atop X_i(T)=1,\,i\in C} F(T)\geq
\sum_{|S|=k,\,X_i(S)=0,i\in A,\atop X_i(S)=1,i\in C} F(S).
\end{equation}
Noting that the sums above have the same number of summands, namely
$$\binom{2k-2|C|+1}{k-|C|+1}\quad\text{and}\quad\binom{2k-2|C|+1}{k-| C|},$$
respectively, we see that \eqref{cna6} is a consequence of the monotonicity  
of $F$.
\end{proof}

Recall from Remark~\ref{r-cn-hn} that CNA+ $\Rightarrow$ Rayleigh/h-NLC+ and 
CNA $\Rightarrow$ h-NLC. By combining Proposition~\ref{RAS} and 
Corollary~\ref{cor-cna+} on the one hand, and Proposition~\ref{pro-h-NLC} and 
Theorem~\ref{th-a-ex} on the other hand, we deduce the following result.

\begin{corollary}\label{RandC}
If $\mu$ is an almost symmetric measure (i.e., its
generating polynomial is symmetric in all but possibly one variable,
cf.~Definition~\ref{d-exc}), then 
\begin{enumerate}
\item $\mu$ is Rayleigh/h-NLC+ $\Longleftrightarrow$ $\mu$ is CNA+;
\item $\mu$ is h-NLC $\Longleftrightarrow$ $\mu$ is CNA.
\end{enumerate}
\end{corollary} 

This shows that Pemantle's Conjecture~\ref{con-cna-R} 
is true for almost 
symmetric measures (recall that by the results in \S \ref{s-sna}, 
Conjecture~\ref{con-cna-R} also holds for strongly Rayleigh measures and
PHR measures).


\section{Negative Results on Negative Dependence}\label{s-countex}

We will now use the results from the previous 
sections, in particular those of \S \ref{Almost}, to construct
counterexamples to the general cases of Pemantle's and Wagner's conjectures
on ULC rank sequences (Conjectures~\ref{c-p}--\ref{cna-ulc}) and   
Pemantle's conjecture
on stochastic domination and truncations 
(Conjecture~\ref{st}). We also answer in the negative  
some related problems as well as Problem~\ref{prob-phr-cna} and 
\cite[Problem 6]{BBCV}.
 
\begin{proof}[Counterexample 1]
Suppose that $\mu$ is a probability measure on $2^{[n+1]}$ whose 
generating polynomial 
$$
f(z_1,\ldots, z_{n+1})= 
z_{n+1}\sum_{k=0}^{n}a_ke_k(z_1,\ldots, z_n)+
\sum_{k=0}^{n}b_ke_k(z_1,\ldots, z_n)
$$
is symmetric in its first $n$ variables. Clearly, 
$$
\De(f)(t)=\sum_{k=0}^{n+1}\binom {n+1} k c_kt^k
=\sum_{k=0}^{n+1}\!\left[\binom n {k-1} a_{k-1}+\binom n k b_k\right]\!t^k, 
$$
where
$$
c_k = \frac{ka_{k-1}+(n+1-k)b_k}{n+1}.
$$
In light of Proposition \ref{RAS} and Corollary \ref{RandC},  
Conjectures~\ref{big} and \ref{cna-ulc} reduce to the following problem in 
the almost exchangeable case.  

\begin{problem}\label{p1}
Let $\{a_k\}_{k=0}^m$, $\{b_k\}_{k=0}^m$ be two non-negative sequences  
satisfying 
\begin{enumerate}
\item $\{\theta a_k+(1-\theta) b_k\}_{k=0}^m$ is LC 
for all $0\leq \theta \leq 1$, and 
\item $a_kb_{k+1} \geq a_{k+1}b_k$ for $0 \leq k \leq m-1$.
\end{enumerate}
Is it true that the sequence 
$$
\left\{ka_{k-1}+(m+1-k)b_k\right\}_{k=0}^{m+1},\,\text{ where }
a_{-1}=b_{m+1}=0, 
$$
is LC? 
\end{problem}

A related problem is the following.

\begin{problem}\label{p2}
Let $\{a_k\}_{k=0}^n$, $\{b_k\}_{k=0}^n$ be two non-negative sequences  
satisfying 
\begin{enumerate}
\item $\{\theta a_k+(1-\theta) b_k\}_{k=0}^n$ is LC 
for all $0\leq \theta \leq 1$, and 
\item $a_kb_{k+1} \geq a_{k+1}b_k$ for $0 \leq k \leq n-1$.
\end{enumerate}
Is it true that the sequence 
$$
\left\{ka_{k-1}+b_k\right\}_{k=0}^{n+1},\,\text{ where }\, a_{-1}=b_{n+1}=0, 
$$
is LC?
\end{problem}

Note that a counterexample $\{a_k\}_{k=0}^n$, $\{b_k\}_{k=0}^n$ to 
the latter problem would actually give a counterexample  to 
the former as well. This is
 because for some large integer $m\geq n$  the sequence $\{d_k\}_{k=0}^m$ 
defined by 
\begin{equation}\label{count-m}
d_k= ka_{k-1}+\frac {m+1-k} m b_k,\,\text{ where }\,     
d_k=0 \mbox{ if } k >n+1,
\end{equation}
will fail to be LC and the sequences  $\{a_k\}_{k=0}^m$ and 
$\{b_k/m\}_{k=0}^m$ 
will therefore constitute a counterexample to Problem~\ref{p1}. 

It remains to find negative solutions to Problems~\ref{p1} and \ref{p2}.   
Set 
$$a_0=2t^2,\quad a_1=2t,\quad a_2=\frac{4}{9}-t,\quad a_3=\frac{2}{3},\quad 
a_4=1,\quad a_5=\frac{2}{3}t,$$
$$b_0=9t^3,\quad b_1=9t^2,\quad b_2=3t,\quad b_3=1,\quad b_4=3,\quad b_5=9-t.$$
For $t$ sufficiently small and positive this is seen to be a negative 
solution to Problem~\ref{p2}. Indeed, 
\begin{equation*}
\begin{split}
&\lambda a_0+b_0=t^2(2\lambda+9t),\,\,\, \lambda a_1+b_1=t(2\lambda+9t),
\,\,\, \lambda a_2+b_2
=\left(\frac{4}{9}-t\right)\lambda+3t,\\
&\lambda a_3+b_3=\frac{2}{3}\lambda +1,\,\,\, \lambda a_4+b_4=\lambda+3,
\,\,\, \lambda a_5+b_5
=\frac{2}{3}t\lambda+9-t,
\end{split}
\end{equation*}
hence
$$(\lambda a_1+b_1)^2-(\lambda a_0+b_0)(\lambda a_2+b_2)
=t^2(2\lambda+9t)\left[\left(\frac{14}{9}+
t\right)\lambda+6t\right]\ge 0,$$
$$(\lambda a_2+b_2)^2-(\lambda a_1+b_1)(\lambda a_3+b_3)
=\lambda^2\left[\left(\frac{4}{9}-t\right)^2
-\frac{4}{9}t\right]+2t\lambda \left(\frac{1}{3}-6t\right)\ge 0,$$
$$(\lambda a_3+b_3)^2-(\lambda a_2+b_2)(\lambda a_4+b_4)
=t\lambda^2+(1-9t)\ge 0,$$
$$(\lambda a_4+b_4)^2-(\lambda a_3+b_3)(\lambda a_5+b_5)
=\left(1-\frac{4}{9}t\right)\lambda^2
+t\ge 0$$
for any $\lambda \ge 0$ if e.g.~$0<t\le \frac{1}{20}$. Moreover,
$$a_0b_1-a_1b_0=0,\quad a_1b_2-a_2b_1=t^2(2+9t)\ge 0,\quad a_2b_3-a_3b_2=
\frac{4}{9}-3t\ge 0,$$
$$a_3b_4-a_4b_3=1,\quad a_4b_5-a_5b_4=9-3t\ge 0$$
for $0<t\le\frac{4}{27}$. However, 
$$
(b_4+4a_3)^2-(b_3+3a_2)(b_5+5a_4)= -\frac{5}{9}+\frac{133}{3}t-3t^2<0 
$$
for $0<t \leq \frac{1}{80}$. Thus, if $0<t\le \frac{1}{80}$ then both 
sequences $\{a_k\}_{k=0}^5$ and $\{b_k\}_{k=0}^5$ are positive  and satisfy 
conditions (1) and (2) in 
Problems~\ref{p1} and \ref{p2}. By taking e.g.~$t=10^{-4}$ one can then 
check that the corresponding sequences $\{a_k\}_{k=0}^m$ and 
$\{b_k/m\}_{k=0}^m$ constructed in \eqref{count-m} produce a 
counterexample to Problem~\ref{p1} for any integer $m\ge 19$. Therefore, 
both Conjecture~\ref{big} and Conjecture~\ref{cna-ulc} fail whenever the 
(total) number 
of variables satisfies $n+1\ge 20$. 

We note that our construction disproves slightly more than 
Conjecture~\ref{big} and Conjecture~\ref{cna-ulc}. Indeed, 
the rank sequence of the measures obtained in this fashion 
(i.e., measures  
on $2^{[m+1]}$ with a generating polynomial of the form 
$$
z_{m+1}\sum_{k=0}^{m}a_ke_k(z_1,\ldots,z_m)+
\sum_{k=0}^{m}\frac{b_k}{m}e_k(z_1,\ldots,z_m),
$$
where $m\ge 19$ and $\{a_k\}_{k=0}^{m}$, $\{b_k\}_{k=0}^{m}$ are 
the sequences constructed above with $t=10^{-4}$)  
is not even SLC, which is weaker than ULC (Definition~\ref{d-ulc}). 

In \cite[Problem 6]{BBCV} it was 
asked whether Conjecture~\ref{big} would hold under the additional assumption
that the Rayleigh measure $\mu$ is decreasing on $2^{[n+1]}$, i.e., 
the coefficients of its generating polynomial $\sum a_Sz^S$ satisfy the
monotonicity condition $a_S\ge a_T$ whenever $S\subseteq T\subseteq [n+1]$. 
This fails as well, since if $f$ is as in the above construction with
$n+1\ge 20$, then the polynomial $f(\lambda z_1,\ldots,\lambda z_{n+1})$ is 
still Rayleigh and  
satisfies the additional monotonicity assumption stated above for all 
sufficiently small $\lambda\in\RR_+$, but its
rank sequence $\{\lambda^{k}r_k\}_{k=0}^{n+1}$ is not ULC.
\end{proof}

\begin{remark}\label{r-other-w}
In \cite{W} it was shown that Conjectures 3.11 and 3.4 in 
{\em loc.~cit.~}are equivalent and would follow from either Conjecture 3.13 or
Conjecture 3.14 in that same paper. Counterexample 1 above shows that all these
conjectures fail.
\end{remark}

Next we show that the above construction actually yields a counterexample to 
Conjecture~\ref{st} as well. 

\begin{proof}[Counterexample 2]
Let $$f(z_1,\ldots,z_{n+1})
= z_{n+1}\sum_{j=0}^n a_je_j(z_1,\ldots, z_n) + 
\sum_{j=0}^n b_je_j(z_1,\ldots, z_n)$$ 
be the $n+1=20$ variable polynomial in the previous counterexample 
and let $\mu$ be the 
corresponding probability measure in $\fP_{n+1}$, so that $\mu$ is CNA+. 
For $\AS=\{ S  \in 2^{[n+1]} : n+1 \in S \}$ we have 
$$
\mu_k(\AS) = \frac {a_{k-1} \binom n {k-1}}{a_{k-1}\binom n {k-1}
+b_k\binom n {k}}, \quad \quad 
\mu_{k+1}(\AS) = \frac {a_{k} \binom n {k}}{a_{k}\binom n {k}
+b_{k+1}\binom n {k+1}}.
$$
Hence, 
$$\mu_k(\AS) \leq \mu_{k+1}(\AS) \quad \Longleftrightarrow \quad a_{k-1} 
\binom n {k-1}b_{k+1}\binom n {k+1} \leq a_{k} \binom n {k}b_k\binom n {k} .$$
This fails for $k=3$ and $0<t<\frac{1}{18}$. 
\end{proof}

As we noted in \S \ref{s24}, a positive answer to Problem~\ref{prob-phr-cna}
would imply Conjecture~\ref{con-cna-R}. However, the answer to 
Problem~\ref{prob-phr-cna} is negative, which we will now prove. 

\begin{proof}[Counterexample 3]
We claim that if $\mu\in\fP_m$ satisfies 
\begin{enumerate}
\item $\mu$ is Rayleigh, 
\item $\mu=\theta \mu_k +(1-\theta)\mu_{k+1}$ for some $0<\theta<1$, 
\item $\mu_k \not \pc \mu_{k+1}$, 
\end{enumerate}
where $\mu_k,\mu_{k+1}$ are the truncations of $\mu$ as in 
Definition~\ref{d-trunc}, then $\mu$ is not PHR (projection of a 
Rayleigh measure with homogeneous generating polynomial, 
cf.~Definition~\ref{phr}). 
To prove this claim, let $g_k$ and $g_{k+1}$ be the generating 
polynomials for $\mu_k$ and $\mu_{k+1}$, respectively, and suppose that 
$$
f(z_1,\ldots,z_\ell)=h_{1}(z_{m+1},\ldots, z_\ell)g_k(z_1,\ldots,z_m)
+h_2(z_{m+1},\ldots, z_\ell)g_{k+1}(z_1,\ldots,z_m),
$$
is a homogeneous Rayleigh polynomial that projects to 
$\theta g_k +(1-\theta)g_{k+1}$ and has a minimal number of 
variables $\ell\ge m+1$. If the same variable, say $z_j$, appears in both 
$h_{1}(z_{m+1},\ldots, z_\ell)$ and $h_2(z_{m+1},\ldots, z_\ell)$, 
then we can differentiate $f$ with respect to $z_j$ and get a new 
homogeneous polynomial in fewer variables that still projects to 
$\theta g_k +(1-\theta)g_{k+1}$. Moreover, if either 
$h_{1}(z_{m+1},\ldots, z_\ell)$ or  
$h_2(z_{m+1},\ldots, z_\ell)$ is the positive sum of several monomials we 
could set a variable equal to zero and again get a polynomial in fewer 
variables with the desired properties. Hence, both  
$h_{1}(z_{m+1},\ldots, z_\ell)$ and $h_2(z_{m+1},\ldots, z_\ell)$ are 
monomials and by the minimality of $\ell$ we have that $\ell=m+1$ and
$$
f(z_1,\ldots,z_{m+1})= \theta z_{m+1}g_k(z_1,\ldots,z_m)
+(1-\theta)g_{k+1}(z_1,\ldots,z_m). 
$$
Since $f$ is a homogeneous Rayleigh polynomial, by Theorem \ref{PHR-NA} 
the corresponding measure $\mu_f\in\fP_{m+1}$ is (strongly) negatively 
associated. As in the proof of Proposition~\ref{stdom} we then deduce that 
$\mu_k  \pc \mu_{k+1}$, which contradicts assumption (3) above. 

It remains to find a measure $\mu$ satisfying properties (1)--(3) above.
If we let $\mu$ have generating polynomial a constant multiple of 
$$
z_{n+1}(a_2e_2(z_1,\ldots, z_n) + a_3e_3(z_1,\ldots, z_n))+ 
b_3e_3(z_1,\ldots, z_n) + b_4e_4(z_1,\ldots, z_n),
$$
where $m=n+1$ and $a_2,a_3, b_3, b_4$ are as in Counterexample 1, we see that 
$\mu$ satisfies all these three properties. 
Clearly, this counterexample fails to have a ULC rank sequence. Therefore, 
we conclude that not even the PHR property implies ULC. 
\end{proof}

\begin{remark}\label{r-strong-NA}
We note that to disprove just the strongest version of 
Pemantle's conjecture on ULC rank sequences (Conjecture~\ref{c-p}) we can 
simply invoke Example~\ref{ex2} in \S \ref{s21}, which combined with
Theorem~\ref{pro-w} shows that
Conjecture~\ref{c-p} actually fails even for symmetric measures. To give yet 
another
counterexample to Conjecture~\ref{c-p} we can 
also employ a construction used in \cite{L2}. Indeed, the measure in  
$\fP_4$ considered in \cite[Theorem 3.5]{L2} is seen to be NA but not ULC. 

Note though that both examples mentioned 
above fail to be NLC (hence also CNA+, CNA, h-NLC+, h-NLC), so they do 
not provide counterexamples to the weaker conjectures (Conjectures~\ref{big}
and \ref{cna-ulc}). 
\end{remark}

\medskip

\noindent
{\bf Acknowledgements.} We would like to thank David Wagner for 
discussions pertaining to this project and the anonymous referee for useful 
comments and suggestions. The first author is 
partially supported by Swedish Research 
Council Grant 621-2005-870. The third author is supported by 
NSF Grant DMS-0301795.


\begin{thebibliography}{99}

\bibitem{ASW} 
M.~Aissen, I.~J.~Schoenberg, A.~Whitney, {\em On the generating functions of 
totally positive sequences I}, J. Analyse Math. {\bf 2} (1952), 93--103.

\bibitem{An}
E.~D.~Andjel, {\em A correlation inequality for the symmetric exclusion 
process}, Ann. Probab. {\bf 16} (1988), 717--721.

\bibitem{ABG}  
M.~F.~Atiyah, R.~Bott, L.~G\.arding, {\em Lacunas for
hyperbolic differential operators with constant coefficients I},
Acta Math. {\bf 124} (1970), 109--189.

\bibitem{BGLS}
H.~H.~Bauschke, O.~G\"uler, A.~S.~Lewis, H.~S.~Sendov, {\em Hyperbolic 
polynomials and convex analysis}, Canad. J. Math. {\bf 53} (2001), 470--488.

\bibitem{BKPV}
J.~Ben Hough, M.~Krishnapur, Y.~Peres, B.~Vir\'ag, {\em Determinantal
Processes and Independence}, Probab. Surv. {\bf 3} (2006), 206--229.

\bibitem{Bor}
J.~Borcea, {\em Spectral order and isotonic differential operators of
Laguerre-P\'olya  type},  Ark. Mat. {\bf 44} (2006), 211--240.

\bibitem{BBS}
J.~Borcea, P.~Br\"and\'en, {\em P\'olya-Schur master theorems for 
circular domains and their boundaries}, to appear in Ann. of Math., 
http://annals.math.princeton.edu.

\bibitem{BBS1}
J.~Borcea, P.~Br\"and\'en, {\em Applications of stable 
polynomials to mixed determinants: Johnson's conjectures, unimodality, and 
symmetrized Fischer products}, Duke Math. J. {\bf 143} (2008), 205--223.

\bibitem{BBS2}
J.~Borcea, P.~Br\"and\'en, {\em Multivariate P\'olya-Schur classification problems in 
the Weyl algebra}, arXiv:math/0606360.

\bibitem{BB}
J.~Borcea, P.~Br\"and\'en, {\em The Lee-Yang and P\'olya-Schur programs I. 
Linear operators preserving stability}, preprint (2008).

\bibitem{BB-2}
J.~Borcea, P.~Br\"and\'en, {\em The Lee-Yang and P\'olya-Schur programs II. 
Theory of stable polynomials and applications}, preprint (2008).

\bibitem{BBCV}
J.~Borcea, P.~Br\"and\'en, G.~Csordas, V.~Vinnikov, {\em P\'olya-Schur-Lax 
problems: hyperbolicity and stability preservers}, 
http://www.aimath.org/pastworkshops/polyaschurlax.html.

\bibitem{BOO}
A.~Borodin, A.~Okounkov, G.~Olshanski, {\em Asymptotics of Plancherel
measures for symmetric groups}, J. Amer. Math. Soc. {\bf 13} (2000), 481--515.

\bibitem{BKKKL}
J.~Bourgain, J.~Kahn, G.~Kalai, Y.~Katznelson, N.~Linial, {\em The influence
of variables in product spaces}, Israel J. Math. {\bf 77} (1992), 55--64.

\bibitem{PB} 
P.~Br\"and\'en, {\em Polynomials with the half-plane property and matroid 
theory}, Adv. Math. {\bf 216} (2007), 302--320.

\bibitem{BP}
R.~Burton, R.~Pemantle, {\em Local characteristics, entropy and limit
theorems for spanning trees and domino tilings via
transfer-impedances}, Ann. Probab. {\bf 21} (1993), 1329--1371.

\bibitem{Ca}
D.~Carlson, {\em Weakly sign-symmetric matrices and some determinantal 
inequalities}, Colloq. Math. {\bf 17} (1967), 123--129.

\bibitem{Chaiken}
S.~Chaiken, {\em A combinatorial proof of the all minors matrix tree theorem}, 
SIAM J. Alg. Disc. Meth. {\bf 3} (1982), 319--329.

\bibitem{COSW} 
Y.~Choe, J.~Oxley, A.~Sokal, D.~G.~Wagner, {\em Homogeneous multivariate 
polynomials with the half-plane property}, 
Adv. Appl. Math. {\bf 32} (2004), 88--187.

\bibitem{CW} 
Y.~Choe, D.~G.~Wagner, {\em Rayleigh Matroids}, Combin. Prob. Comput. {\bf 15} 
(2006), 765--781.

\bibitem{Con}
J.~B.~Conrey, {\em The Riemann Hypothesis}, Notices Amer. Math. Soc. {\bf 50}
(2003), 341--353.

\bibitem{CC}
T.~Craven, G.~Csordas, {\em Jensen polynomials and the T\'uran and Laguerre
inequalities}, Pacific J. Math {\bf 136} (1989), 241--260.

\bibitem{DVJ}
D.~J.~Daley, D.~Vere-Jones, An Introduction to the Theory of Point Processes.
Springer, New York, 1988.

\bibitem{DJR}
D.~Dubhashi, J.~Jonasson, D.~Ranjan, {\em Positive influence and negative 
dependence}, Combin. Probab. Comput. {\bf 16} (2007), 29--41. 

\bibitem{DPR}
D.~Dubhashi, V.~Priebe, D.~Ranjan, 
{\em Negative Dependence Through the FKG Inequality}, Technical Report 
RS-96-27, BRICS Report Series, Basic Research In Computer Science, \AA rhus, 
Denmark, 1996; webversion available at http://citeseer.ist.psu.edu/352490.html.

\bibitem{DR}
D.~Dubhashi, D.~Ranjan, 
{\em Balls and bins: A study in negative dependence}, Random Struct. 
Algorithms {\bf 13} (1998), 99--124.

\bibitem{EK}
S.~N.~Ethier, T.~G.~Kurtz, Markov Processes: Characterization and
Convergence. Wiley, 1986.

\bibitem{FJ}
S.~M.~Fallat, C.~R.~Johnson, {\em Determinantal Inequalities: Ancient History
and Recent Advances}, Contemp. Math. {\bf 259} (2000), 199--212.

\bibitem{feder}
T.~Feder, M.~Mihail, {\em Balanced matroids}, in ``Proceedings of the 24th 
Annual ACM (STOC)'', ACM Press, New York, 1992.

\bibitem{FKG}
C.~M.~Fortuin, P.~W.~Kasteleyn, J.~Ginibre, {\em Correlation inequalities on 
some partially ordered sets}, Commun. Math. Phys. {\bf 22} (1971), 89--103.

\bibitem{GK}
F.~R.~Gantmacher, M.~G.~Krein, {\em Oscillation matrices and kernels and small
vibrations of mechanical systems}, Gostechizdat, 1950.

\bibitem{grace} 
J.~H.~Grace, {\em The zeros of a polynomial}, Proc. Cambridge Philos. Soc. 
{\bf 11} (1902), 352--357.

\bibitem{grimmett} 
G.~Grimmett, The Random Cluster Model. Springer, 2006.  

\bibitem{Gu1}
L.~Gurvits, {\em  A proof of hyperbolic van der Waerden conjecture: 
the right generalization is the ultimate simplification}, arXiv:math/0504397.

\bibitem{Gu}
L.~Gurvits, {\em Hyperbolic polynomials approach to van der 
Waerden/Schrijver-Valiant like Conjectures: sharper bounds, simpler proofs 
and algorithmic
applications}, to appear in Proc. STOC 2006,  arXiv:math/0510452.


\bibitem{guler} 
O.~G\"uler, {\em Hyperbolic polynomials and interior point 
methods for convex programming}, Math. Oper. Res. {\bf 22} (1997), 350--377.

\bibitem{garding} 
L.~G\aa rding, {\em An inequality for hyperbolic polynomials}, 
J. Math. Mech. {\bf 8} (1959), 957--965. 

\bibitem{HLP} 
G.~H.~Hardy, J.~E.~Littlewood, G.~P\'olya, Inequalities. 2nd ed., Cambridge 
Univ. Press, Cambridge, UK, 1988.

\bibitem{HM}
O.~Holtz, {\em $M$-matrices satisfy Newton's inequalities}, Proc. Amer. Math. 
Soc. {\bf 133} (2005), 711--717.

\bibitem{HS}
O.~Holtz, H.~Schneider, {\em Open problems on GKK $\tau$-matrices}, Linear 
Algebra Appl. {\bf 345} (2002), 263--267. 

\bibitem{HSt}
O.~Holtz, B.~Sturmfels, {\em Hyperdeterminantal relations among symmetric
principal minors}, J. Algebra {\bf 316} (2007), 634--648.

\bibitem{Hag}
O.~H\"aggstr\"om, {\em Random-cluster measures and uniform spanning trees},
Stoch. Proc. Appl.  {\bf 59} (1995), 267--275.

\bibitem{hormander}
L. H\"ormander, Notions of Convexity. Progr. Math. {\bf 127}. 
Birkh\"auser, Boston, MA, 1994.

\bibitem{JJP}
G.~James, C.~Johnson, S.~Pierce, {\em Generalized matrix function 
inequalities on $M$-matrices}, J. London Math. Soc (2) {\bf 57} (1998), 
562--582.

\bibitem{J-DP}
K.~Joag-Dev, F.~Proschan, {\em Negative association of random variables with
applications}, Ann. Stat. {\bf 11} (1983), 286--295.

\bibitem{Joh}
K.~Johansson, {\em Discrete orthogonal polynomial ensembles and the Plancherel
measure}, Ann. of Math. (2) {\bf 153} (2001), 259--296.

\bibitem{Joh2}
K.~Johansson, {\em Determinantal processes with number variance saturation},
Commun. Math. Phys. {\bf 252} (2004), 111-148.

\bibitem{KN}
J.~Kahn, M.~Neiman, {\em Negative correlation and log-concavity}, 
arXiv:math/0712.3507.

\bibitem{Kar}
S.~Karlin, Total Positivity. Vol. I. Stanford Univ. Press, Stanford, CA, 
1968.




\bibitem{Le}
B.~Ja.~Levin, Distribution of Zeros of Entire Functions. 
Transl. Math. Monogr. {\bf 5}, Amer. Math. Soc., Providence, R.I., 1980. 

\bibitem{LS}
E.~H.~Lieb, A.~D.~Sokal, {\em A General Lee-Yang Theorem for One-Component and 
Multicomponent Ferromagnets}, Commun. Math. Phys. {\bf 80} (1981), 153--179.

\bibitem{L0}
T.~M.~Liggett, Interacting Particle Systems. Springer, 1985.

\bibitem{L1}
T.~M.~Liggett, {\em Ultra log-concave sequences and negative dependence}, 
J. Combin. Theory Ser. A {\bf 79} (1997), 315--325.

\bibitem{L00}
T.~M.~Liggett, Stochastic Interacting Systems: Contact, Voter and
Exclusion Processes. Springer, 1999.

\bibitem{L2}
T.~M.~Liggett, {\em  Negative correlations and particle systems}, Markov 
Proc. Rel. Fields {\bf 8} (2002), 547--564.

\bibitem{L3} 
T.~M.~Liggett, {\em Distributional limits for the symmetric 
exclusion process}, to appear in Stoch. Proc. Appl, arXiv:math/0710.3606. 

\bibitem{Lyons}
R.~Lyons, {\em  Determinantal probability measures}, Publ. Math. Inst. 
Hautes \'Etudes Sci. {\bf 98} (2003), 167--212.

\bibitem{Lyons-Peres}
R.~Lyons, Y.~Peres, Probability on Trees and Networks. Book in progress, 
webversion available at 
http://mypage.iu.edu/\symbol{'176}rdlyons/prbtree/prbtree.html.

\bibitem{Lyons-Steif}
R.~Lyons, J.~Steif, {\em Stationary determinantal processes: phase 
multiplicity, Bernoullicity, entropy, and domination}, Duke Math. J. {\bf 120}
(2003), 515--575.

\bibitem{Mark} 
K.~Markstr\"om, {\em Negative association does not imply log-concavity of the 
rank sequence}, J. Appl. Probab. {\bf 44} (2007), 1119--1121; updated version:
\texttt{http://abel.math.umu.se/\~{}klasm/}.

\bibitem{MO}
A.~Marshall, I.~Olkin, Inequalities: Theory of Majorization and Its 
Applications. Academic Press, New York, 1979.

\bibitem{Mason} 
J.~H.~Mason, {\em Matroids: unimodal conjectures and Motzkin's theorem,} 
Combinatorics (Proc. Conf. Combinatorial Math., Math. Inst., Oxford, 1972), 
pp. 207--220. Inst. Math. Appl., Southend-on-Sea, 1972.

\bibitem{New}
C.~Newman, {\em Normal fluctuations and the FKG inequalities}, Commun. Math. 
Phys. {\bf 74} (1980), 119--128.

\bibitem{OR}
A.~Okounkov, N.~Reshetikhin, {\em Correlation function of Schur process with
applications to local geometry of a random 3-dimensional Young diagram},
J. Amer. Math. Soc. {\bf 16} (2003), 581--603.

\bibitem{P}
R.~Pemantle, {\em Towards a theory of negative dependence}, J. Math. Phys. 
{\bf 41} (2000), 1371--1390.

\bibitem{RS}
Q.~I.~Rahman, G.~Schmeisser, Analytic Theory of Polynomials. London
Math. Soc. Monogr. (N.~S.) {\bf 26}, Oxford Univ.~Press, New York,
NY, 2002.

\bibitem{rota}
G.-C.~Rota, D.~Sharp, {\em Mathematics, Philosophy, and Artificial 
Intelligence: a Dialogue with Gian-Carlo Rota and David Sharp}, 
Los Alamos Science, Spring/Summer 1985.

\bibitem{semple-welsh}
C.~Semple, D.~J.~A.~Welsh, {\em Negative correlation in graphs and matroids}, Combin. Prob. Comput. {\bf 15} (2006), 765--781.

\bibitem{seymour-welsh}
P.~D.~Seymour, D.~J.~A.~Welsh, {\em Combinatorial applications of an inequality from statistical mechanics}, Math. Proc. Camb. Phil. Soc. {\bf 77} (1975), 485--495.

\bibitem{sok}
A.~D.~Sokal, {\em The multivariate Tutte polynomial (alias Potts model)
for graphs and matroids}, in ``Surveys in Combinatorics, 2005'' 
(B.~S.~Webb, ed.), Cambridge Univ. Press, Cambridge, UK, 2005. 

\bibitem{szego}
G.~Szeg\"o, {\em Bemerkungen zu einem Satz von J.~H.~Grace \"uber die Wurzeln 
algebraischer Gleichungen}, Math. Z. {\bf 13} (1922), 28--55.

\bibitem{W}
D.~G.~Wagner, {\em Negatively correlated random variables and Mason's 
conjecture for independent sets in matroids}, to appear in Ann. Combin., 
arXiv:math/0602648.

\bibitem{W2}
D.~G.~Wagner, {\em Matroid inequalities from electrical network theory},
Electron. J. Combin. {\bf 11} (2005), A1 (17pp).

\bibitem{walsh}
J.~L.~Walsh, {\em On the location of the roots of certain types of 
polynomials}, Trans. Amer. Math. Soc. {\bf 24} (1922), 163--180.

\end{thebibliography}
\end{document}